\documentclass{article}
\usepackage[utf8]{inputenc}
\usepackage{amsmath,amssymb,dsfont,amsfonts}
\usepackage{mathrsfs}
\usepackage{graphics}
\usepackage{graphicx}
\usepackage{color}
\newcommand{\Prob}{\mathbb P}
\newcommand{\R}{\mathbb{R}}

\newcommand{\N}{\mathbb N}

\newcommand{\proof}{{\noindent \bf Proof:} }
\newcommand{\eop }{ \hfill $\Box$ }

\newtheorem{example}{Example}

\newtheorem{theorem}{Theorem}[section]
\newtheorem{proposition}[theorem]{Proposition}
\newtheorem{corollary}[theorem]{Corollary}
\newtheorem{remark}[theorem]{Remark}
\newtheorem{definition}[theorem]{Definition}

\usepackage[bottom=2cm,top=2.5cm,left=2.5cm,right=2.5cm]{geometry}

\title{Decomposition of discontinuous flows of diffeomorphisms: jumpings, geometrical and topological aspects}
\author{Lourival Lima\footnote{Mathematics Department, Escuela Superior Politécnica del Litoral, ESPOL. Vía Perimetral 5, Guayaquil, Ecuador. E-mail: lourodri@espol.edu.ec. 
Research supported by ESPOL.} \ \ and\ \ Paulo Ruffino\footnote{Mathematics Department, State University of Campinas, UNICAMP. Rua Sérgio Buarque de Holanda 651. 13083-659 Campinas, SP, Brazil. E-mail: ruffino@unicamp.br. 
Research partially supported by CNPq 305212/2019-2, FAPESP 2020/04426-6.} }

\date{}

\begin{document}

\maketitle

\begin{abstract} Let $M$ be a compact manifold equipped with a pair of complementary
foliations, say horizontal $\mathcal{H}$ and vertical $\mathcal{V}$. In Melo, Morgado and Ruffino (Disc Cont Dyn Syst B,  2016, 21(9)) it is proved that if a semimartingale $X_t$ has a finite number of jumps in compact intervals then, up
to a stopping time $\tau$, a stochastic flow of local diffeomorphisms in  $M$ driven by $X_t$  can be decomposed into a process  in the Lie group of diffeomorphisms which fix the leaves of $\mathcal{H}$  composed with a process in the Lie group of 
diffeomorphisms which fix the leaves of $\mathcal{V}$. Dynamics at the discontinuities of $X_t$ here are interpreted in the Marcus sense as in Kurtz, Pardoux and Protter \cite{KPP}. 
Here we enlarge the scope of this geometric decomposition and consider flows driven by arbitrary semimartingales with jumps and show explicit equations for each component. Our technique is based in an extension of the Itô-Ventzel-Kunita formula for stochastic flows with jumps. Geometrical and others topological obstructions for the decomposition are also considered: e.g. an index of attainability is introduced to measure the complexity of the dynamics with respect to the pair of foliations. 
\end{abstract}

\section{Introduction}

In this article we study a class of geometrical decomposition of flows of diffeomorphisms generated by some dynamical systems which covers most of the discontinuous case in the recent literature. As a first model of this geometrical decomposition consider a dynamical system (flow of local diffeomorphisms) $\varphi_t$ on an $n$-dimensional differentiable manifold $M$ endowed with a pair of foliations $\mathcal{H}$ and $\mathcal{V}$, say. The notation here stands for a visualization as \textit{horizontal} and \textit{vertical} foliations.  Assume that this pair of foliations is  \textit{complementary} in the sense that the sum of their dimensions equals $n$, and at each point $x\in M$, the direct sum of the tangent spaces of the leaves passing through $x$, $L^1_x $ and $L^2_x $ in  $\mathcal{H}$ and $\mathcal{V}$ respectively, is the tangent space $T_xM$. The first model of decomposition we are interested here is a local description of the behaviour of open sets in a neighbourhood of an initial condition $x_0\in M$ in terms of the leaves of $\mathcal{H}$ and $\mathcal{V}$: We study the possibility of writing $\varphi_t = \eta_t \circ \psi_t$, where $\eta_t$ are diffeomorphisms which fix the leaves of $\mathcal{H}$ and  $\psi_t$ are diffeomorphisms which fix the leaves of $\mathcal{V}$. Hence, if $x\in M$ is in the domain, then the trajectory  $t \mapsto \eta_t (x)$ lies on the leave $L^1_ x$ and $t \mapsto \psi_t (x)$ lies on the leave $L^2_ x$. Up to changing of coordinates, this decomposition corresponds to generalize the classical Cartesian decomposition from single trajectories, to the dynamics of open sets. In the next sections this simple model is going to be made precise and will be extended to further analytical and geometrical context.

\bigskip

 We establish geometrical and analytical conditions for the existence of such decomposition of flows along the leaves of this bifoliated space. Some of these conditions can be intrinsically related to the manifold, some examples are given in Sections 3 and 4, where we define the index of attainability, which reflects the level of topological obstruction for the existence of the decomposition.  We also state a technique to perform a decomposition of the form: 
\begin{eqnarray}
\varphi_{t}(x_0) = \big( \eta^k_{t \vee s_k}\circ\psi^k_{t\vee s_k} \big) \circ ... \circ \big( \eta^2_{s_2}\circ \psi^2_{s_2} \big) \circ \big( \eta^1_{s_1}\circ \psi^1_{s_1}\big) (x_0) ,
\label{fc1}
\end{eqnarray}
where $\eta^i$ and $\psi^i$, $i=1,2 \ldots, k$ are purely vertical and horizontal components respectively. We call this factorization as an \textbf{alternate decomposition}.

\bigskip

The article is organized in the following way: In the next Section we present formally the definition of biregular atlas in a manifold with two complementary foliations. We also recall previous results on decomposition for continuous flows, in particular the cascade decomposition: when one has more than two complementary foliations. In Section 3  we prove a generalized It\^o-Ventzel-Kunita formula in the context of the Marcus approach of the jumps cf. Kurtz, Pardoux and Protter \cite{KPP}. This formula is crucial on the composition of distinct flows generated by semimartingales with jumps. Here, in particular we extend the previous result for finite jumps in compact interval, \cite{Melo2}, to an arbitrary number of jumps; it turns then possible to apply the formula for a larger variety of noises which includes, for example, Lévy noise, see e.g. Applebaum \cite{D.Applebaum}, Protter \cite{Protter}, Oksendal and Sulem \cite{Oksendal}, among others.  Related results with a different technique is also proposed for L\'evy process, among others, in Hartmann and Pavlyukevich \cite{Hartmann-Pavlyukevich}, Qiao and Duan \cite{Qiao-Duan}, Qiao \cite{Qiao}. Different from the approach of these previous articles, here our formula focus on the dynamics of composition of flows associated to  general semimartingales with jumps. In this Section we prove the decomposition Theorem \ref{disc version} for semimartingales with jumps. In fact, the main decomposition theorem is proved in the more general context of complementary distributions, recall that if these distributions are integrable then they generate complementary foliations. In Section 4 we propose a method of \textit{alternate decomposition} to overcome the problem of explosion time in the single step (horizontal-vertical) decomposition. The \textit{index of attainability} of the system is introduced: it indicates the topological complexity of the trajectories with respect to the complementary foliations (Def. \ref{def: attainability} and Proposition \ref{Prop: Attan}). In Section 5 we apply our technique in principal fibre bundles over homogeneous spaces. Here the horizontal distribution is given by a connection $\omega$ in the Lie group which acts transitively on homogeneous space. 

\section{Foliations and previous results on continuous flows}

Given an $n$-dimensional smooth manifold $M$, a foliation $\mathcal{F}$ of dimension $1 \leq k < n$ in $M$ is a partition of $M$ into immersed connected submanifolds of dimension $k$, called the leaves of $\mathcal{F}$ with local foliated chart. That is, locally  $(M, \mathcal{F})$ is diffeomorphic to open sets of $\mathbb{R}^n = \mathbb{R}^k \times \mathbb{R}^{n-k}$, in such a way that the leaves have constant second coordinate. In fact, a foliation $(M, \mathcal{F})$ is identified with a foliated atlas which is coherent along the leaves in the following sense:

\begin{definition} \emph{Let $M$ be a smooth $n$-dimensional manifold. A (smooth) $k$-dimensional foliated atlas $\mathcal{A}$ of $M$ is a maximal atlas  on $M$ which satisfies:}
\begin{itemize} 
    \item[1)] If $(U,\alpha) \in \mathcal{A}$, then $\alpha(U) = U_1 \times U_2 \subset \mathbb{R}^k \times \mathbb{R}^{n-k}$ for $U_1$, $U_2$ open subsets of $\mathbb{R}^k$ and $\mathbb{R}^{n-k}$ respectively. 
    
    \item[2)]  Given two local charts $(U,\alpha), (V,\beta)\in \mathcal{A}$, with $U \cap V \neq \emptyset$, then the change of coordinate map is given by $\alpha \circ \beta^{-1}(x,y) = (h_1(x,y),h_2(y))$, for some smooth maps $h_1$ and $h_2$ in the appropriate domain.
\end{itemize} 
\end{definition}
A foliated atlas $\mathcal{A}$  is said to be \emph{regular} if it is locally finite and for any foliated chart $(U, \alpha)\in \mathcal{A}$, the closure of its domain $\bar{U}$ is a compact set contained in $V$, the domain of another foliated chart $(V, \beta)$. The sets $\alpha^{-1} (B, \{y\})\subset M$, for $(U,\alpha)\in \mathcal{A}, B\subset \R^k$ open set such that $ (B, \{y\})\subset \alpha(U)$ are called \emph{plaques} of the atlas.

\bigskip

Consider the  equivalence relation in $M$ given by $x\sim y$ if and only if there exists a finite sequence of plaques $P_0, P_1, \ldots, P_p$ with $x\in P_0$, $y\in P_p$ and $P_i\cap P_{i-1}\neq \emptyset$ for all $i=1, \ldots , p$. The equivalent classes here determine a one-to-one correspondence between regular foliated atlases and the leaves $F$ of a foliated manifold (see e.g. \cite[Thm.  1.2.18]{Candel}).
 Given a point $p\in M$, the unique leaf of the foliation  passing through  $p$ is denoted by $L_p$.  The set $\displaystyle{\mathcal{F}(S)= \cup_{p \in S}L_p}$ is called the  \textit{saturation of $S$} by $\mathcal{F}$.


\begin{definition}
\emph{We say that a maximal atlas of a foliated space $(M, \mathcal{F})$ is transversely orientable  if the cotangent bundle of the leaves in the tangente bundle of $TM$ is orientable. Precisely: for all change of coordinates $\phi_1 \circ \phi_2^{-1} : \phi_2(U_1 \cap U_2) \longrightarrow \mathbb{R}$ we have that
\begin{eqnarray}
\det \frac{\partial (y_{n-k+1}, ..., y_n)}{\partial (x_{n-k+1}, ..., x_n)} > 0.
\end{eqnarray}
where  $ (y_1, ..., y_n)=\phi_1 \circ \phi_2^{-1} (x_1, \ldots, x_n)$.}
\end{definition}

\begin{definition} \label{Def: biregular charts}
\emph{Consider a bifoliated manifold $M$, i.e. it is endowed with a pair of complementary foliation in $M$, say horizontal $\mathcal{H}$ and vertical $\mathcal{V}$ foliation. An atlas  $\mathcal{A}$ is  biregular on $(M,\mathcal{H}, \mathcal{V})$ if $\mathcal{A}$ is foliated and regular for $\mathcal{H}$ and $\mathcal{V}$ simultaneously, moreover, given two coordinate systems $(U,\alpha_1)$ and $(V,\alpha_2)$, with $U \cap V \neq \emptyset$, the change of coordinate map is given by $\alpha_1 \circ \alpha_2^{-1}(x,y) = (h_1(x),h_2(y))$, for some smooth maps $h_1$ and $h_2$ in the appropriate domain.}
\end{definition}

The existence of a biregular atlas is straightforward, see e.g. \cite[Lemma 5.1.4]{Candel}. Given an initial condition $x_0 \in M$. Unless otherwise stated, we are going to assume that the horizontal foliation $(M, \mathcal 
 {H})$ is tranversely orientable. This is not quite a restriction since 
if $\mathcal{H}$  is not transversely orientable, it can be lifted to a transversely-oriented foliation on a double covering of $M$, see e.g. \cite[Prop. 3.5.1]{Candel}.

\begin{remark}[Decomposition of a fixed diffeomorphism:]
\label{characterization of decomposition}
\emph{
A paradigmatic example of the decomposition we are treating here can be described for a single diffeomorphism, say $\phi:U\subset \R^n \rightarrow \phi(U) $ in the following way: for $1\leq k\leq n$ consider the factorization of the state space $\R^n = \R^k \times \R^{n-k}$. Consider the Cartesian pair of foliations where: the horizontal foliation $\mathcal{H}$ of $\R^n$ is given by horizontal affine subspaces  $\R^k \times \{y\} $, for $y\in \R^{n-k}$; analogously, the vertical foliation $\mathcal{V}$ of $\R^n$ is given by the vertical affine subspaces  $ \{x\} \times \R^{n-k} $, for $x\in \R^{k}$. In these coordinates, the group of diffeomorphisms which fix the horizontal leaves are given by (up to the domain) $\psi(x,y)= (\psi^1(x,y), y)$ for a differentiable map $\psi^1:V \rightarrow \R^k$. Analogously, the group of diffeomorphisms which fix the vertical leaves are given by $\eta(x,y)= (x, \eta^2(x,y))$ for a differentiable map $\eta^2:V \rightarrow \R^{n-k}$.  In these coordinates, the original diffeomorphism $\phi$ can be written as  $\displaystyle{\phi = \left(\phi^1(x,y), \phi^2(x,y) \right)}$, with $(x,y) \in \R^{n-k} \times \R^k$. Given a point $x_0 \in U$, it follows directly by the inverse function theorem that there exists a unique (reducing the domain if necessary) decomposition $\displaystyle{\phi = \eta \circ \psi}$ in a neighbourhood of $x_0$, where $\eta$ and $\psi$ are horizontal and vertical preserving diffeomorphisms, if and only if 
\begin{eqnarray}
\label{prop1}
\det \frac{\partial \phi^2(x_0)}{\partial y} \neq 0.
\end{eqnarray}}

\bigskip

\emph{Applying this characterization for a flow of diffeomorphism $\varphi_t$, one can guarantee the local existence of decomposition $\displaystyle{\varphi_t = \eta_t \circ \psi_t}$ up to a stopping 
 time $\tau$, where $\eta_t$ and $\psi_t$ are families of diffeomorphisms preserving horizontal and vertical components respectively where
\begin{eqnarray}
\tau = \sup \left\{t > 0; \det\frac{\partial \phi^2_s(x,y)}{\partial y} \neq 0 \hspace{0.2cm} \mbox{for all} \hspace{0.2cm} 0 \leq s \leq t \right\}.
\end{eqnarray}
For more details see Melo et al \cite{Melo1}. In some cases, there exists some degree of compatibility of the vector fields with the complementary distributions in such a way that the decomposition presented above holds for all time $t$, i.e. the determinant of equation (\ref{prop1}) never vanishes, see \cite{Melo2}. A basic but important example is a linear system on $\mathbb{R}^n \setminus \{0\}$ endowed with spherical and radial foliations, note that in this specific case, the system sends radial leaves into radial leaves, therefore the decomposition holds for all time. Another example in this context is the derivative flow $\displaystyle{\varphi_{t*}: T_{x_0}M \longrightarrow  T_{\varphi_t(x_0)}M}$, for $x_0 \in M$, in the linear frame bundle $\displaystyle{\pi: BM \longrightarrow M}$. Note that $\varphi_{t*}$ is an isomorphism between the fibres $\pi^{-1}(x_0)$ and $\pi^{-1}(\varphi_t(x_0))$, therefore $\varphi_{t*}$ has a decomposition $\varphi_{t*} = \eta_t \circ \psi_t$ for all time $t\geq 0$, where $\psi_t$ is linear in the tangent space $T_{x_0}M$ into itself and $\eta_t$ is horizontal with respect to a connection in $TM$.}
\end{remark}
\eop

\bigskip

The following basic result is crucial on determining the topology of the attainable  sets which we are going to introduce in Section 4. It is a nontrivial result if one considers e.g. non-compact or dense leaves in a compact foliated space. We adapt its proof from \cite{Camacho} into our context. 

\begin{theorem}[Bi-foliated uniform transversality]
\label{Thm: Uniform transversality}

Consider $(M, \mathcal{H}, \mathcal{V})$, a manifold with complementary foliations $\mathcal{H}$ and $\mathcal{V}$. Fix a leaf $F$, say, in $\mathcal{H}$. Given two points $p,q\in F$, let $V_p,  V_q \in \mathcal{V}$ be the vertical leaves passing thorough $p$ and $q$ respectively. Then, there exist open sets in the intrinsic topology $D_1 \subset V_p$, $D_2 \subset V_q$ with $p \in  D_1$, $q \in D_2$ and a diffeomorphism $f: D_1 \rightarrow D_2$ such that $f(L\cap D_1)= L \cap D_2$ for every horizontal leaf $L$ in $\mathcal{H}$.
\end{theorem}

\proof Consider local biregular charts $\varphi_p: U_p \rightarrow U_1 \times U_2$  and 
$\varphi_q: U_q \rightarrow \tilde{U}_1 \times \tilde{U}_2$ in a neighbourhood of $p$ and $q$ respectively, with $U_1, \tilde{U}_1 \subset \R^k$, $U_2, \tilde{U}_2 \subset \R^{n-k}$ and $\varphi_p(p)=\varphi_q(q)=(0,0)$. By the uniform transversality theorem, see e.g. \cite[Thm. 3, Ch.III]{Camacho} there exist submanifolds $N_1$ and $N_2$, with $p\in N_1$ and $q\in N_2$ transverse to $F$ and a diffeomorphism $\tilde{f}: N_1 \rightarrow N_2$ such that $\tilde{f}(L\cap N_1)=L\cap N_2$
for all horizontal leaf $L$ .

We have to show that $N_1$ and $N_2$ above can be chosen as open sets $D_1$ and $D_2$ of the vertical leaves $V_p$ and $V_q$. Since $N_1$ is transverse to $F$ at $p$, then the derivative at $p$ of the non-linear projection  $\psi_p:=\varphi_p^{-1} (\{0\} \times \pi_2 \circ \varphi_p\circ i): N_1 \rightarrow V_p $ is an isomorphism between the tangent spaces $T_p N_1$ and $T_p V_p$, where $i:N_1 \rightarrow M$ is the inclusion and $\pi_2: \R^n \rightarrow \R^{n-k}$ is the projection. By the classical local inverse theorem, there exists an open set $\tilde {D_1}$ where the restriction of $\psi_p $ is a diffeomorphism. By the same argument, we have that there exists an open set $\tilde{D}_2$ where the restriction of $\psi_q:=\varphi_q^{-1} (\{0\} \times \pi_2 \circ \varphi_q \circ i) $ is also a diffeomorphism. 

The diffeomorphism $f: D_1 \rightarrow D_ 2$ of the statement is given by $f=\psi_q  \circ \tilde{f} \circ \psi_p^{-1}$ with $D_1 = \psi_p (\tilde{D_1}\cap \tilde{f}^{-1}(\tilde{D_2}))$ and $D_2 = \psi_q ( \tilde{f}(\tilde{D_1})\cap \tilde{D_2})$.

\eop

\subsection{Previous results on decomposition of continuous flows}

Many studies focusing on distinct geometrical decomposition of flows have examined dynamics driven by Stratonovich equations, see e.g. \cite{Catuogno}, \cite{Melo2}, \cite{Liao1} among many others. The technique in these articles employs intrinsic calculus of continuous semimartingales on manifolds. Here we recall the main results of existence of the decomposition of $\varphi_t$, a continuous stochastic flow of (local) diffeomorphisms generated by an autonomous stochastic Stratonovich equation on $M$. 


\bigskip

More general than foliations, suppose that $M$ is endowed  with a pair of  differentiable complementary distributions. That is, we have two complementary  differentiable assignment of subspaces in each fibre of $TM$. More precisely, locally we have:  $\Delta^H: U \subset M \rightarrow Gr_{k}(M)$ and $\Delta^V: U \subset M \rightarrow Gr_{(n-k)}(M)$ respectively, where $U \subset M$ is a connected open set and $\displaystyle{Gr_k}(M) = \bigcup_{x \in M}Gr_k(T_xM)$ is the Grasmannian bundle. They are complementary in the sense that  $\displaystyle{\Delta^H(x) \oplus \Delta^V(x) = T_xM}$, for all $x \in U$. Consider the following Stratonovich stochastic differential equation on the manifold $M$:

\begin{eqnarray}
\label{eq6}
dx_t = \sum_{i=0}^kX_i(x_t)\circ dW^i_t,
\end{eqnarray}
with initial condition $x_0 \in M$, where $X_0, X_1, \ldots , X_k$ are smooth with bounded derivative vector fields on $M$, $(W^1_t, \ldots, W^k_t)$ is a Brownian motion on $\mathbb{R}^k$, and $(W^0_t) = t$.  We suppose that this probabilistic structure is well defined over an appropriate filtered probability space $(\Omega, \mathcal{F}, (\mathcal{F})_{t\geq 0}, \mathbb{P})$. Let $\varphi_t: \Omega \times M \rightarrow M$ be the stochastic flow associated to the diffusion generated by equation (\ref{eq6}). If we assume that the derivatives of the vector fields are bounded, then $\varphi_t$ exists for all $t \geq 0$. We denote by $\exp \{tX\} \in \mbox{Diff}(M)$, the local flow of diffeomorphisms associated with $X$, and by $\mbox{Diff}(\Delta , M)$, the group of diffeomorphisms generated by exponentials of vector fields in $\Delta$, where $\Delta$ is a distribution in $M$. In other words:  
\begin{equation*}
\mbox{Diff}(\Delta, M) = \mbox{cl} \left\lbrace \exp\{t_1X_1\} \ldots \exp\{t_nX_n\}, \mbox{ with } X_i\in\Delta, t_i\in\mathbb{R}, \forall n\in\mathbb{N} \right\rbrace .
\end{equation*}

\noindent The Lie group generated by all smooth vector fields $\mbox{Diff}(TM, M)$ contains two important Lie subgroups: the subgroup generated by all purely horizontal vector fields, denoted by $\mbox{Diff}(\Delta^H, M)$ and the subgroup generated by all purely vertical vector fields  denoted by $\mbox{Diff}(\Delta^V, M)$. If the distributions are integrable, hence generate foliations, then the intersection of these subgroups is the identity and each element of these groups preserves the leave of the corresponding foliation.


\begin{definition}
\emph{A complementary pair of distributions $\Delta^H$ and $\Delta^V$ in $M$ is said to preserve transversality along the orbits of the subgroup $\emph{Diff}(\Delta^H, M)$ (acting on $TM$), if for any element $F \in \mathrm{Diff}(\Delta^H, M)$ we have  that  $DF \Delta^V\left(F^{-1}(x)\right)\cap \Delta^H(x) = \{0\}$.   }
\end{definition}

\begin{theorem}[Decomposition of continuous flows] Let $\Delta^{H}$ and $\Delta^{V}$ be two complementary distributions on the manifold $M$ which preserves transversality along $\emph{Diff}(\Delta^H, M)$. If $\varphi_t$ is the flow of diffeomorphism generated by equation (\ref{eq6}) then there exists, up to a stopping time, a decomposition:
\begin{eqnarray}
\label{comp1}
\varphi_t = \eta_t \circ \psi_t.
\end{eqnarray}
Where $\eta_t$ is a diffusion (solution of an autonomous SDE) in $\emph{Diff}(\Delta^H, M)$ and $\psi_t$ is a process (solution of a non-autonomous SDE) in $\emph{Diff}(\Delta^V, M)$.
\label{thm1}
\end{theorem}

For a proof, see \cite[Thm 2.2]{Catuogno}.  An interesting fact is that if $\Delta^{H}$ and $\Delta^{V}$ are both involutive distributions, then locally, the decomposition (\ref{comp1}) is unique. In fact, in this case, in a neighbourhood $U_{x} \subset M$ the holonomy of the foliations vanishes, in other words  $\text{Diff}(\Delta^H, U_{x}) \cap  \text{Diff}(\Delta^V, U_{x}) = I$. 

\bigskip

Exploring this decomposition with dimensions or codimensions of the foliations increasing one by one, also motivated by the possibility of studying flow of diffeomorphisms over sections of  flag bundles, we have a \textbf{cascade decomposition}. Let $(\Delta^{H}_1, \ldots, \Delta^{H}_k)$ and $(\Delta^{V}_k, \ldots, \Delta^{V}_1)$ be two smooth sections of a flag bundle over $M$. i.e. we have that $\Delta^{H}_1 \subset \ldots \subset \Delta^{H}_k$ and $\Delta^{V}_1 \supset \ldots \supset \Delta^{V}_k$. It is worth mentioning that for $\text{dim}M = n$, the sequence $(\Delta^{H}_1, \ldots, \Delta^{H}_k)$ is a section of the maximal flag manifold if $\text{dim}\Delta^{H}_{i+1} - \text{dim}\Delta^{H}_i = 1$ and $k=n$. Suppose that $\Delta^{H}_i$ and $\Delta^{V}_i$ are complementary in each point of the manifold $M$, then, since those distributions are actually sections of a flag bundle, it yields that $\text{diff}(\Delta^{H}_i, M) \subseteq \text{diff}(\Delta^{H}_{i+1}, M)$ and $\text{diff}(\Delta^{V}_{i+1}, M) \subseteq \text{diff}(\Delta^{V}_{i}, M)$ for each $i = 1, \ldots, k-1$. Applying Theorem \ref{thm1} recurrently in this flag structure, we have the following

\begin{corollary} \label{Cor: decomp flag}  Given the complementary flag structure as describe above, locally the original flow $\varphi_t$ of equation (\ref{eq6})  can be decomposed up to an explosion time as:

\begin{equation}
    \label{Eq: decomp para flag}
\varphi_t = \eta^1_t \circ \eta^2_t \circ \ldots \eta^k_t \circ \psi_t,
\end{equation} 
with the following property: the composition of the first $i$-th components $(\eta^1_t \circ \ldots \eta^i_t)$ solves an autonomous SDE in  $  \emph{diff}(\Delta^H_i, M)$ $i=1, \ldots k$ and the composition of the last $(k-i+1)$ components $(\eta^{i+1}_t \circ \ldots \eta^k_t \circ \psi_t)$ solve a (non autonomous) SDE in $  \emph{diff}(\Delta^V_{i}, M)$. This decomposition is unique when the distributions of the flag bundles are involutive.
\end{corollary}

\begin{example}\emph{ Consider a local coordinate system  $\phi=(\phi_1, \ldots,
\phi_n): U \subset M \rightarrow \R^n$ over $M$. Then $\phi$ establishes in $U$ a biregular atlas of 1-dimensional complementary foliations such that the corresponding distributions satisfy the geometrical assumptions  of the previous corollary. Hence, given a flow of (local)
diffeomorphisms $\varphi_t$, up to an explosion time, locally there exists the decomposition
\[
 \varphi_t=\eta^1_t \circ  \eta_t^2\ldots \circ \eta_t^n
\]
where each diffeomorphism $\eta^i_t$
preserves the $j$-th coordinates for all $j\neq i$. Note that in this example, the component $\psi_t$ of Corollary \ref{Cor: decomp flag} degenerates to the identity. For details of the proof, see \cite{Catuogno}}.

\end{example}

\section{A generalization of Itô-Ventzel-Kunita formula}

Many interesting flows of diffeomorphisms are not continuous, for example, those generated by semimartingales with jumps, including Lévy noise and many others, see e.g. Kurtz, Pardoux and Protter \cite{KPP}. The main result of this section, Theorem \ref{teo: IVK} is a generalization of the Itô-Ventzel-Kunita formula for  composition of discontinuous flows  generated by noise with jumps in the sense of  Marcus equation as in  \cite{KPP}. Here, we enlarge the scope of the previous result in \cite{Melo2}, allowing infinite many jumps in compact intervals.

\subsection{Stratonovich SDE with jumps on manifolds (SDEJ)}
\noindent Initially we recall the definition and main properties of an Stratonovich SDE with jumps (SDEJ) in the sense of Marcus equation in an Euclidean space, among others see e.g. \cite{KPP}. We consider a standard probability structure given by $(\Omega, \mathcal{F}, \mathcal{F}_t, \Prob)$ a complete filtered probability space, where our semimartingales are defined. Let $Z = \{ Z_t, t\geq 0\}$ be a $k$-dimensional semimartingale, with $Z_0 = 0$, and let  $[Z,Z] = [Z^j,Z^m]$ be the quadratic variation matrix which can be decomposed into $[Z,Z] = [Z,Z]^c + [Z,Z]^d$, where $[Z,Z]^c$ and $[Z,Z]^d$ represent the continuous and purely discontinuous parts  respectively. The process $\Delta Z_t= Z_{t^+}-Z_{t^-}$ is the pure jump component of the c\`adl\`ag semimartingale $Z_t$. Let $X \in C^1(\mathbb{R}^n; \mathcal{L}(\mathbb{R}^k, \R^n))$ with $X=(X^1, \ldots , X^n)$ being $k$ smooth with bounded derivative vector fields in $\R^d$, i.e. $X^i=Xe_i$, for $e_i$ elements of the canonical basis in $\R^n$. Given an initial condition  $x_0\in \R^n $, an SDEJ is given by

\begin{eqnarray}
\label{Eq. KPP}
x_t = x_0 + \int_0^tX(x_s) \diamond dZ_s,
\end{eqnarray}
where, shortly speaking: the continuous part of the solution corresponds to the classical Stratonovich equations and the jump part are performed along fictitious time, i.e. jumps of time $\Delta Z$ along the deterministic flow $\phi(X, x, t)$ of $X$ or equivalently, jumps of time one of $X \Delta Z_t$. Equation (\ref{Eq. KPP})
has a unique solution up to a stopping time $\tau$. Precisely the solution $x_t$ satisfies:
\begin{eqnarray} \label{Eq: sol Marcus}
x_t &=& x_0 + \int_0^t X(x_{s^-}) \, dZ_s + \frac{1}{2}\int_0^t X'X(x_s d [Z,Z]_s^c) \\
 & & \\
 & & \sum_{0<s\leq t }\left\{ \phi(X\Delta Z_s, x_{s^-}) - x_{s^-} - X(x_{s^-}) \Delta Z_s \right\}.
\end{eqnarray}
For the convergence of the series on the right hand side, further properties and details of this formula, see  \cite[eq. (2.2)]{KPP}.  

\bigskip

The integration performed in equation (\ref{Eq. KPP}) is not properly an integral since it does not apply to a class of integrand, it is only defined for $X(x_t)$ in the context of the corresponding SDEJ. In this sense, to show a change of variable formula for equation (\ref{Eq. KPP}), one has to extend the meaning of this integral. More precisely given a differentiable function $g: \R^n \rightarrow \R^k$, to overcome this problem, Kurtz, Pardoux and Protter introduce an extended definition of Marcus equation at another category of integrand (tensor here) along the trajectories of (\ref{Eq. KPP}) to allow the following formula, \cite[Def. 4.1]{KPP}:
\[
\int_0^t  g(x_s) \diamond dZ_s := \int_0^t g(x_{s^-}) dZ_s  + \frac{1}{2}\, \mathrm{Tr}\int_0^t g'(x_s)  (dZ, dZ)^c  (x_{s}) X^t
\]
\begin{equation} \label{Def: Marcus-Ito}
\hspace{5cm} + \sum_{0<s\leq t}\left[ \int_0^1 
\left( g(\phi(X\Delta Z_s,x_{s^-}, u)) - g(x_{s^-}) \right) du  \right] \Delta Z.
\end{equation}
Note that if $g$ is the vector fields $X$ itself, the definition of equation (\ref{Def: Marcus-Ito}) above coincides with the original definition of Marcus solution: in fact, in this case $\phi(X\Delta Z_s, x_ {s^-}, u)$ itself is a primitive of the first summand $g(\phi(X\Delta Z_s,x_{s^-}, u)) \Delta Z_s$ of the integral in the second line. 

With this extended interpretation above we have the following change of variable formula:
\begin{proposition}
\label{Prop: mudança de variaveis}
If $x_t$ is the solution of equation ($\ref{Eq. KPP}$), then for a diffeomorphism $f \in C^2(\mathbb{R}^n)$:
\begin{eqnarray} \label{Eq: Ito for Marcus}
f(x_t) = f(x_0) + \int_0^t f'(x_s)X(x_s) \diamond dZ_s,  \ \ t \geq 0.
\end{eqnarray}
\end{proposition}
For a proof, see \cite[Prop. 4.2]{KPP}. Observe that since conjugacy of (deterministic) flow is equivalent of the conjugacy of the corresponding vector fields then, in the formula (\ref{Def: Marcus-Ito}) we have that  $f \psi(X\Delta Z_s, x_{s^-}, \cdot)= \psi(f'X \Delta 
Z_s, f(x_{s^-}), \cdot )$. Hence $f(x_t)$ is solution of the Marcus equation
\[
y_t = y_0 + \int_0^t f' X (y_s) \diamond Z_t.
\]
I.e., solutions $x_t$ and $y_t$ are conjugated by the diffeomorphism $f$.
\bigskip

Marcus equation can be considered in a differentiable manifold $M$: in fact, in \cite{KPP} they use an embedding argument of $M$ in an Euclidean space to define it. A natural intrinsic extension can be given using the change of variable formula of Proposition \ref{Prop: mudança de variaveis} above: let $X\in \mathcal{C}^1(M; \mathcal{L}(\mathbb{R}^k, TM))$, such that for each $x \in M$ the linear map $X(x)$ sends a vector $z\in\mathbb{R}^k$ into $X(x)z\in T_xM$. Assume that the vector field $X$ is smooth with bounded derivative on $M$ and consider the equation
\begin{equation}
dx_t = X(x_t)\diamond dZ_t, \quad x(0)=x_0.
\label{Eq de Marcus}
\end{equation}
We say that  $x_t \in M$, $t \geq 0$ is a solution of equation (\ref{Eq de Marcus}) if for any $f \in C^2(M)$,
\[
f(x_t)= f(x_0) + \int_ 0^t f'(x_s) X(x_ s) \ \diamond dZ_s
\]
in the sense of equation (\ref{Def: Marcus-Ito}), i.e. it is interpreted as:   
\begin{eqnarray}
\label{interpretacao}
f(x_t) &=& \int_0^t df X(x_{s^-})\, dZ_s + \frac{1}{2}\, \mathrm{Tr}\int_0^t \nabla^2f (XdZ, XdZ)^c (x_{s}) \nonumber \\
&+& \sum_{0<s\leq t}\left[f(\phi(X\Delta Z_s,x_{s^-}, 1)) - f(x_{s^-}) - df X (x_{s^-})\right] \Delta Z_s.\hspace*{5cm}
\end{eqnarray}
The first term on the right hand side of equation (\ref{interpretacao}) is a standard It{\^o} integral of the predictable process $df X^j(x_{s^-})$ with respect to the semimartingale $Z_t$. The second term is the Stieltjes integral of the Levi-Civita connection applied in the derivative of the function $f$, with respect to the continuous part of the quadratic variation of $Z_t$.
In the third term: $\phi(X\Delta Z_s, x_{s^-}, 1)$ indicates the solution
at a fictitious time $t=1$ of the ODE generated by the vector field $X\Delta Z_s$ and initial condition $x_{s^-}$ (if last entry `1' is ommited, it means we are considering time one). Thus, the jumps of this equation occurs in deterministic directions. Usual regularities conditions over the linear map $X(x)$ and its derivatives guarantee the existence of a unique Stratonovich flow of diffeomorphisms $\varphi_t$, which is solution of equation (\ref{Eq de Marcus}). Moreover, for an embedded submanifold $M$ in an Euclidean space, a support theorem \cite[Prop. 4.3]{KPP} states that the solution still remains on the manifold after a jump.  

\bigskip

Again, using the It\^o formula (\ref{Eq: Ito for Marcus}), in  local charts of $M$ the equation (\ref{Eq de Marcus}) has the following local solution: 

\begin{proposition}
Let  $\alpha: U \subset M \longrightarrow \mathbb{R}^n$ be a local coordinate system of $M$ in a neighbourhood of $x_0$. Let $Y=\alpha_* (X)$ be the induced vector fields in $\R^n$ and consider the SDEJ in $\R^n$ given by
\begin{equation} \label{Eq: conjugacy}
dy_t = Y(y_t)\diamond dZ_t, \quad y(0)=\alpha (x_0).
\end{equation}
Then, up to a stopping time, the solution flows of  equation (\ref{Eq de Marcus}) and equation (\ref{Eq: conjugacy}) are conjugate by $\alpha$, i.e. $x_t= \alpha^{-1}(y_t)$.

\end{proposition}
\proof
It follows straightforward by the It\^o formula of Proposition (\ref{Prop: mudança de variaveis}) up to a stopping time. We check that the solution $x_t$ on the manifold is independent of the coordinate system: In fact, given another local chart $\beta$ with appropriate domain (in space and time)  the diffeomorphism $\beta \circ \alpha^{-1}$ conjugates the Marcus equations generated by $\alpha_*X$ and $\beta_*X$, by formula (\ref{Eq: Ito for Marcus}).

\eop

\subsection{Itô-Ventzel-Kunita for Stratonovich SDEJs}

In order to prove the main theorem of this Section (Theorem \ref{disc version}) we need  an extension of the classical It\^{o}-Ventzel-Kunita formula (infinite jump version). As we said in the paragraph before Proposition \ref{Prop: mudança de variaveis}, the Marcus integral only makes sense in a SDEJ. So, for the same reason that in \cite{KPP} they have to define equation (\ref{Def: Marcus-Ito}) for the composition with differentiable mapping, here we have to define what is the integral of a vector field induced by the derivative of a flow:  Precisely, let $X$ and $Y$ be two smooth vector fields on $\mathbb{R}^n$ and consider $F_t$ and $G_t$, flows of diffeomorphisms generated by the SDEJs $dF_t = X(F_t) \diamond dZ_t$ and $dG_t = Y(G_t) \diamond dZ_t$, with respect to the same semimartingale $Z_t$. Cf. \cite{Melo2}, we define the following integral: 
\begin{eqnarray}
\nonumber \int_0^t(F'_{s}Y(G_s))\diamond dZ_s &:=& \int_0^t(F'_{s^-}Y(G_{s^-}))dZ_s \\
\nonumber &+& \frac{1}{2}\int_0^t(X'(Y(G_s))+F_{s*}(Y'Y))d[Z,Z]_s^c \\
\nonumber &+& \sum_{0 \leq s \leq t}\Big\{\phi(X\Delta Z_s,F_{s^-}(\phi(Y\Delta Z_s, G_{s^-})) - \phi(X\Delta Z_s,F_{s^-}(G_{s^-})) \\
&-&  (F'_{s^-}Y(G_{s-})) \Delta Z_s \Big\}, \label{Def: gen IKV}
\end{eqnarray}
where the first term on the right hand side is the Itô integral of $F_{s*}Y(G_{s^-})$ with respect to $Z_s$. In the second term, note that $(X'(Y(G_s))+ F_{s*}(Y'Y)) = d F_{s*}Y(G_s)(F_{s*}Y(G_s))$, so the second integral corresponds to the finite variation such that its continuous part satisfies the classical Itô-Ventzel-Kunita. On the last summation, the expression $\phi(X\Delta Z_s,F_{s^-}(\phi(Y\Delta Z_s, G_{s^-}))$ has the following geometrical meaning: at the jump time $s \in [0,t]$,  flow $G_{(\cdot)}$ jumps before the jump of the flow $F_{(\cdot)}$. Note that if $F_t=Id$ then the definition of the solution of Marcus equation (\ref{Eq: sol Marcus}) is recovered.

\bigskip

The summation in the definition given by expression (\ref{Def: gen IKV}) is absolutely convergent, in fact, applying Taylor's formula in the map   $$\displaystyle{u \longrightarrow \phi(X\Delta Z_s,F_{s^-}(\phi(Y\Delta Z_s, G_{s^-},u),1),}$$ around $u=0$
one finds
\begin{eqnarray*}
\phi(X\Delta Z_s,F_{s^-}(\phi(Y\Delta Z_s, G_{s^-},1),1) &=& \phi(X \Delta_sZ,F_{s^-}(G_{s^-})) + (F_{s^*}Y(G_s)) \Delta Z_s \\
&& \\
&+& \frac{1}{2} \frac{\partial ^2 \phi}{\partial  u^2 }(\theta_1,\theta_2)\Delta Z_s \Delta Z^t_s,
 \end{eqnarray*}
where the Lagrange second order remainder has $\theta_1,\theta_2 \in (0,1)$.
Therefore, 
\begin{eqnarray*}
\sum_{0<s \leq t}\Big|\phi(X\Delta_s Z,F_{s^-}(\phi(Y\Delta Z_s, G_{s^-})) - \phi(X\Delta Z_s,F_{s^-}(G_{s^-})) 
-  (F_{s*}Y(G_s))(F_s \circ G_s) \Delta Z_s \Big| 
\end{eqnarray*}
\vspace{-0.5cm}
\begin{eqnarray*}
&\leq& \sup_{0<s \leq t}\frac{1}{2}\Big|(X'(Y(F_{s^-}))+F_{s*}(Y'Y))f(\phi(X\Delta Z_s,F_{s^-}(\phi(Y\Delta_s Z, G_{s^-},\theta_1),\theta_2))\Big|\sum_{0<s \leq t}|\Delta Z_s|^2  \\
&&\\
&\leq& K \sum_{0<s \leq t}|\Delta Z_s|^2,
\end{eqnarray*}
for a constant $K$. Convergence holds since the sum of squares of the jumps of a general semimartingale is finite a.s..

\bigskip

Next theorem states an extension of It\^{o}-Ventzel-Kunita for general semimartingales. In this context, an infinite number of jumps may occur.

\begin{theorem}[Itô-Ventzel-Kunita for Stratonovich SDEJ]
Suppose that $F_s$ and $G_s$ are solutions of SDEJs driven by a general semimartingale $Z_t$ with respect to smooth with bounded derivative vector fields $X$ and $Y$ on $\mathbb{R}^d$ respectively, then for $t \in [0,T]$,   
\begin{eqnarray}
\label{IVK-formula}
F_s(G_s) = F_0(G_0) + \int_0^t X(F_s(G_s)) \diamond dZ_s + \int_0^t F_{s*}(Y(G_s)) \diamond dZ_s
\end{eqnarray}
\label{teo: IVK}
\end{theorem}

\proof It is well known that formula (\ref{IVK-formula}) holds if $Z_t$ is continuous for each $t \in [0,T]$, see e.g. Kunita \cite[Thm. 8.3]{KPP}. Moreover, in \cite{Melo2} it was proved that if $Z$ has a finite number of jumps in compact intervals, then formula (\ref{IVK-formula}) also holds. Hence, we only have to prove the formula when the semimartingale $Z_t$ jumps infinitely many times in bounded intervals. Compared with the technique used in \cite{Melo2}, in this case the problem arises at times which are accumulation points of jumps. Here, we overcome this problem by splitting the set of jump times of $Z_t$, into two disjoint subsets, $A=A(\epsilon,t)$ and $B=B(\epsilon,t)$, such that $A$ is finite,  $\displaystyle{\sum_{s\in B}(\Delta Z_s)^2\leq \epsilon}$, and $A \cup B$ contains every jump time of $Z_t$. Let $Z^A_t = Z_t- \Delta Z_{1_B(t)}$, i.e., $Z^A_t$ only has discontinuities when  $t\in A$. Note that $Z^A_t$ converges to $Z_t$ uniformly. Consider $F^A_t$ and $G^A_t$, solutions of equations $dF^A_t = X(F^A_t) \diamond dZ^A_t$ and $dG^A_t = Y(G^A_t) \diamond dZ^A_t$, respectively. Hence, formula (\ref{IVK-formula}) holds for all $t\in [0,T]$, i.e. 
\begin{eqnarray}
\label{formula da prova}
F^A_s(G^A_s) = F_0(G_0) + \int_0^t X(F^A_s(G^A_s)) \diamond dZ_s + \int_0^t F^A_{s*}(Y(G^A_s)) \diamond dZ^A_s.
\end{eqnarray}

If $\epsilon$ goes to zero, $A(\epsilon,t)$ expands to sets which are collecting more and more jumps on the interval $[0,t]$ and the set $B$, on the other hand, reduces to a set of arbitrarily small jumps. Precisely,  $F^A_t \longrightarrow F_t$ and $G^A_t \longrightarrow G_t$ uniformly a.s, moreover, by Taylor's formula on the time variable we have that 
\begin{eqnarray*}
\Big|\Big|\int_0^t X(F^A_t(G^A_t))\diamond dZ^A_s - \int_0^tX(F_s(G_s)) \diamond dZ_s    \Big|\Big| &\leq & \sum_{s \in B}\Big|\Big| \Big\{ \phi(X\Delta Z_s, F_{s^-}(G_{s^-})) - F_{s^-}(G_{s^-})\\
&-& X(F_{s^-}(G_{s^-}))\Delta Z_s \Big\} \Big|\Big| \\
&\leq & \sum_{s \in B} K_1 \Big|\Big| \Delta Z_t  \Big|\Big|^2 \leq K_1 \,  \epsilon.
\end{eqnarray*}
For a positive constant $K_1$. Similarly, for the same reason,

\begin{eqnarray*}
\hspace{-3.2cm}\Big|\Big|\int_0^t F^A_{s*}(Y(G^A_s))\diamond dZ^A_s - \int_0^t F_{s*}(Y(G_s)) \diamond dZ_s \Big|\Big|  \leq  K_2 \sum_{s \in B} \Big|\Big| \Delta Z_t  \Big|\Big|^2 \leq K_2 \, \epsilon.
\end{eqnarray*}
Therefore, when $\epsilon$ goes to zero, expression (\ref{formula da prova}) converges uniformly to formula  ($\ref{IVK-formula}$) a.s..

\eop

\bigskip

Next corollary establishes a Leibniz formula for Stratonovich SDEJs. The proof follows directly from Theorem \ref{teo: IVK}. This formula is fundamental  in the next section in order to compute explicit expressions for the components of the decomposition. 

\bigskip

\begin{corollary}[Leibniz formula]
\label{corol1}
Let $F_t$ and $G_t$ be flows generated by Stratonovich SDEJs with respect to a general semimartingale $Z_t$. Then 
\begin{eqnarray}
\diamond d(F\circ G)_t = \diamond d(F_t) G_t + (F_t)_{*}  \diamond dG_t.
\end{eqnarray}
\end{corollary}

By Proposition (\ref{Prop: mudança de variaveis}) and local coordinate arguments we can easily extend all results in this section for a Riemannian manifold.

\subsection{Decomposition of flows of diffeomorphism with jump components}


Now we return to the problem of decomposing a flow of diffeomorphism $\varphi_t \in \mbox{Diff}(M)$, into components in subgroups of $\mbox{Diff}(M)$. In this context, one of the components is again a flow in this space (hence satisfies the Markov property if the original one does).  This kind of decomposition appears in the literature, for example, in Bismut \cite{Bismut}, Kunita \cite{HK}, Ming Liao \cite{Liao}, among others. In this section, we study the existence of this decomposition, when $\varphi_t$ is a solution of a Marcus equation driven by a general semimartingale $Z$ and give equations of its components in each subgroup of $\mbox{Diff}(M)$. In the main result of this section, we prove that the diffeomorphism $\varphi_t$ is decomposable up to a stopping time $\tau$ even if $Z$ perform infinitely many jumps in compact intervals. 

\bigskip

Using Lie group terminology, the dynamics of the stochastic flow $\varphi_t$ which is the (local) solution of the Stratonovich SDEJ (\ref{Eq de Marcus}) lies in 
$\mathrm{Diff}(M)$, the connected Lie group of diffeomorphisms of $M$ generated by vector fields. As any other dynamics generated by vector fields, the flow $\varphi_t$ can be written as the following right invariant SDEJ in $\mathrm{Diff}(M)$:
\begin{eqnarray}
d\varphi_t = R_{\varphi_t^*}X \diamond dZ_t,
\end{eqnarray}
where $R_{\varphi_t^*}$ is the derivative of the right translation in the Lie group $\mbox{Diff}(M)$. 
A first version of the theorem of decomposition for processes with jump components was proved in \cite[Prop. 1]{Melo2} restricted to systems whose driven process $Z_t$ jumps only a finite number of times on compact intervals.  Here, Corollary \ref{corol1} allows us to extend to systems with arbitrary (countable) number of jumps. In section 3.1, we illustrate for linear systems each step in proof of this theorem:

\begin{theorem}[discontinuous version]
\label{disc version}
The stochastic flow of local diffeomorphisms $\varphi_t$ generated by the SDEJ (\ref{Eq de Marcus}) can be 
decomposed locally, up to a
stopping time, as
\[
 \varphi_t = F_t \circ G_t,
\]
where $F_t$ is solution of an (autonomous) SDEJ in $\emph{Diff}(\Delta^H, M)$ and $G_t$ is
a process in $\emph{Diff}(\Delta^V, M)$.
The decomposition is unique if $\Delta^H$ and $\Delta^V$ are integrable.
\end{theorem}
\proof
For each horizontal diffeomorphism $F$, let $\widetilde{X}_F$ be the element in the Lie algebra of the group $ \mbox{Diff}(\Delta^H, M)$, i.e., a horizontal vector field which lives in $\Delta^H$, given by:
\begin{eqnarray}
\label{horizontal vector fields}
\widetilde{X}_F(x) = X(x)-V(x) \ \in \Delta^H,
\end{eqnarray}
where $V(x)$ is the unique vector field in the subspace $dF (\Delta^V(F^{-1} (x))) $ such that $\widetilde{X}_F$ is horizontal in $T_xM$.
The first component $F_t$ of the statement is the solution of the Marcus equation in $\mbox{Diff}(\Delta^H,M)$ given by:
\begin{eqnarray}
\label{jump horizontal}
dF_t =  R_{F_{t*}}\widetilde{X}_{F_t} \diamond dZ^i_t. 
\end{eqnarray}
In fact, since the vector fields are horizontal and are translated by horizontal diffeomorphims, then the solution of this right invariant equation is also horizontal.

\bigskip

For the second component, take  $G_t = F^{-1}_t \circ \varphi_t$. We only have to proof that $G_t$ is vertical. In order to find a Marcus equation whose solution flow is $G_t$, we apply Corollary (\ref{corol1}) (Itô-Ventzel-Kunita for general semimartingales). Hence,
\begin{eqnarray}
\label{second jump component}
dG_t = F^{-1}_{t*} \diamond d\varphi_t + \diamond dF^{-1}_t(\varphi_t).
\end{eqnarray}
On the other hand, note that, by Corollary \ref{corol1} again we have that:
\begin{eqnarray}
\label{inverse of first jump component}
dF^{-1}_t = - L_{F^{-1}_{t*}}  \widetilde{X}_{F_t} \diamond dZ_t. 
\end{eqnarray}
By (\ref{second jump component}) and (\ref{inverse of first jump component}), it follows that:
\begin{eqnarray}
\label{jump vertical}
dG_t = \sum_{i=0}^m \mbox{Ad}(G^{-1}_t)(V_i(G_t))(F_t) \diamond dZ^i_t.
\end{eqnarray}

Finally note that $G_t \in \mathrm{Diff}(\Delta^V, M)$ since $\mbox{Ad}(F^{-1}_t)(V_i)(x)$ are in $\Delta^2(x)$ by construction. If the distributions  $\Delta^H$ and $\Delta^V$ are integrable then they are the tangent bundle of a pair of complementary foliations, in this case, there is no holonomy, that is the intersection $\mathrm{Diff}(\Delta^H, M) \cap \mathrm{Diff}(\Delta^V, M) = \{Id\}$.

\eop

The stopping time stated in the Theorem is not infinite in general (as it will be in the alternate decomposition in Section 4). An illustrative and simple example is the following: consider the pure rotation in $\R^2$ given by $\displaystyle{dx_t = Ax_t \diamond dZ_t}$, where $A$ is the skew-symmetric matrix with entries 1 and $-1$. The flow of this system is given by:
\begin{eqnarray*}
\varphi_t = \left( \begin{array}{cc} \cos Z_t & -\sin Z_t \\ \sin Z_t & \cos Z_t \end{array} \right). 
\end{eqnarray*}
Taking the Cartesian horizontal and vertical foliations, locally (in time) the flow $\varphi_t$ is decomposed by 
\begin{eqnarray}
\label{decomp-jumps}
\hspace{-1cm}\left( \begin{array}{cc} \cos Z_t & -\sin Z_t \\ \sin Z_t & \cos Z_t \end{array} \right) &=& \left( \begin{array}{cc} \sec(Z_t) & -\tan(Z_t ) \\ 0 & 1 \end{array} \right)\left( \begin{array}{cc} 1 & 0 \\ \sin(Z_t) & \cos(Z_t) \end{array} \right),
\end{eqnarray}
which does not exist whenever
$Z_{t} \in \{\frac{\pi}{2} +k\pi, k \in \mathbb{Z}\}$; cf. equation  (\ref{prop1}).

\subsection{Linear Systems} 

In this section we explore further the linear systems in $\R^n$ with the Cartesian foliations and find conditions for existence of decomposition for all time $t\geq 0$. Consider the following  SDEJ:
\begin{equation} \label{Eq: Linear}
  dx_t = A\, x_t \diamond  dZ_t,
\end{equation}
with $x_0\in \R^n $ and, in this context, $Z_t$ a semimartingale on the real line. The Itô-Ventzel-Kunita result in the previous section shows that the  fundamental linear solution flow of (\ref{Eq: Linear}) is the exponential 
\begin{equation} \label{Eq: flow of Linear}
  F_t = \exp{ \{ A Z_t \} }.
 \end{equation}

The decomposition we are interested here is 
\[
F_t = \eta_t \circ \psi_t
\]
such that $\eta_t \in \mbox{Diff}(\Delta^H, M)$ and $\psi_t \in \mbox{Diff}(\Delta^V, M)$. In the general theory of last section, the components $\eta_t$ and $\psi_t$ are not necessarily linear, even in quite symmetric situations: Consider for example, the pair of foliations  in $\R^n \setminus \{0\}$  given by radial and spherical coordinates; in this case neither the radial neither the spherical components are linear. Nevertheless, in the case of the Cartesian pair of foliation $\R^k \times \R^{\ell}$, we do have that $\eta_t$ and $\psi_t$ are actually linear. In fact, writing the fundamental solution in coordinates, we have that 
\[
F_t = \left(  \begin{array}{ll}
\Big( F_1 (t) \Big)_{k\times k}     & \Big( F_2(t)\Big)_{k \times \ell} \\
    & \\
    \Big( F_3 (t)\Big)_{\ell \times k} & \Big( F_4(t)\Big)_{\ell \times \ell}
\end{array}
\right).
\]
Since $\eta_t$ does not change the last $\ell$ coordinates, the diffeomorphisms $\psi_t$ must perform all the transformations on the last $\ell$ coordinates, i.e. it must have  the following form:
\[
\psi_t = \left(  \begin{array}{ll}
   \Big( 1d \Big)_{k\times k}     &  0  \\
   & \\
     F_3(t) & F_4 (t)
\end{array}
\right).
\]
This implies that both $\eta_t$ and $\psi_t$ are linear.

\bigskip





Observe that the vertical component in fact lies in the following Lie group:

\[
\psi_t \in G_V= \left\{  g\in Gl(n, \R); g= \left(\begin{array}{ll}
\Big( 1d \Big)_{k\times k}     &  0  \\
     g_3  & \Big( g_4 \Big)_{\ell \times \ell}
\end{array}
\right)  \right\}
\]
whose Lie algebra is given by the vector space generated by 
\[
\left(\begin{array}{cc}
\Big( 0 \Big)_{k\times k}     & 0 \\
& \\
     \Big(  * \Big)  & \Big( * \Big)_{\ell \times \ell}
\end{array}
\right),
\]
Where $(*)$ denotes arbitrary real entries with the appropriate dimensions. For the horizontal component: 
\[
\eta_t \in G_H= \left\{  g\in Gl(n, \R); g= \left(\begin{array}{ll}
\Big( g_1 \Big)_{k\times k}     &  g_2  \\
     0  & \Big( 1d \Big)_{\ell \times \ell}
\end{array}
\right)  \right\}
\]
with Lie algebra generated by 
\[
\left(\begin{array}{cc}
\Big( * \Big)_{k\times k}     & \Big(  * \Big)_{k \times \ell} \\
& \\
     0  & \Big( 0 \Big)_{\ell \times \ell}
\end{array}
\right).
\]

\bigskip

We illustrate here the calculations of the proof of Theorem \ref{disc version} in order to find the SDEJs for the  submatrices $g_1, g_2$ and  $g_3, g_4$ of $\eta_t$ and $\psi_t$ respectively. Consider $\pi_2: \mathbb{R}^k\times \mathbb{R} ^{\ell}\rightarrow \mathbb{R}^{\ell}$ the projection on the second (vertical) subspace. From the proof of Theorem (\ref{disc version}) we have that
\begin{eqnarray}
\label{caso linear}
V(\eta, \cdot) = \eta \circ \pi_2 \circ A (\cdot).
\end{eqnarray}
In fact, note that $V(\eta, \cdot)$ is in the image of the vertical component by $\eta$ and that  $\pi_2 V(\eta, \cdot) = \pi_2 A (\cdot)$. From (\ref{caso linear}) and   equations (\ref{jump horizontal}) and (\ref{jump vertical}) we find the autonomous equation:
\[
\label{caso linear horizontal}
d\eta_t = (1d - \eta_t \circ \pi_2)A \, \eta_t \ \diamond dz_t,
\]
and the well expected nonautonomous vertical diffeomorphisms:
\[
d\psi_t =   \pi_2\,  A\,  \eta_t \circ \psi_t \ \ \diamond dz_t.
\]

Writing the matrix of coefficients in blocks as
\[
A =: \left(  \begin{array}{ll}
\Big( A_1 \Big)_{k\times k}     & \Big( A_2 \Big) \\
    \Big( A_3 \Big) & \Big( A_4 \Big)_{\ell \times \ell}
\end{array}
\right),
\]
we can calculate explicitly each constituent submatrices of $\eta_t$ and $\psi_t$:
\begin{eqnarray}
d g_1(t) &=& \big[  A_1 \ g_1(t) - g_2(t)\ A_3 \ g_1(t) \big] \ \diamond  dz_t \nonumber\\
d g_2(t) &=& \big[ A_1 g_2(t)+ A_2 - g_2(t)\, A_4 - g_2\, A_3\, g_2(t) \big] \ \diamond dz_t, \nonumber \\
d g_3(t) &=& \big[ A_3\,g_1 + A_3\, g_2\, g_3 + A_4 g_3 \big] \ \diamond dz_t \nonumber\\
d g_4(t) &=& \big[ A_3\, g_2\, g_ 4 + A_4\, g_4 \big] \ \diamond dz_t. 
\label{Eq: constituent}
\end{eqnarray}
 All the terms of equations above are linear except those quadratics  which are all multiplied by the submatrices $A_3$. Hence, if $A_3=0$ then there exists solutions for these equations, i.e. there exists the decomosition for all time $t \geq 0$. Cf. elementary example of rotation after the proof of the Theorem, where $A_3=1$.
 The next Proposition allow us to extend the scope of this decomposition. Before that, we define the following notation: given two complementary subspaces $E_1 \oplus E_2 = \R^n$, we denote by $\mathcal{F}(E_1)$  and $\mathcal{F}(E_2)$ the corresponding pair of complementary affine foliations in $\R^n$ generated by the corresponding affine subspaces.

\begin{proposition}
Consider a linear SDEJ in $\R^n$, with $n>2$, given by 
\begin{equation}\label{Eq: linear}
dx_t = A\, x_t \diamond dz_t.
\end{equation}
Let  $r=\# \{ \mbox{real eigenvalues with multiplicities} \} $ and consider positive integers  $a\leq r$ and $b\leq (n-r)/2$. Then, for any dimension of the form $a+2b$, there exist a pair of affine foliations $\mathcal{F}(E_1)$, $\mathcal{F}(E_2)$ generated by complementary subspaces $E_1$ and $E_2$, $\dim E_1 =a+2b$, such that the decomposition of the flow of equation (\ref{Eq: linear}) exists for all time $t\in [0,T]$. \end{proposition}
\proof
Let $A$ be represented in its canonical real Jordan form as $A = P J P^{-1}$, such that the nilpotent components are above the main diagonal. The coordinate transformation $y=P\, x$ establishes the conjugate Marcus system:
\begin{equation*}
dy_t = J\, y_t \diamond dz_t.
\end{equation*}
If $n>2$, it is possible to express $J$ as 
\[
J = \left(  \begin{array}{ll}
\Big( J_1 \Big)_{k\times k}     & \Big( J_2 \Big) \\
    \Big( J_3 \Big) & \Big( J_4 \Big)_{\ell \times \ell}
\end{array}
\right)
\]
with $k=a + 2b$ and its complementary $\ell=n-k$, such that the submatrix $(J_3)_{\ell \times k} = 0$. Observe that  $a$ represents the number of real eingenvalues in the block $J_3$ and $b$ is the number of pairs of conjugate nonreal eigenvalues in the same block. Hence, equation (\ref{Eq: constituent}) guarantees that there is no explosion in the decomposition of $y_t$. By conjugacy, there is also no explosion in the decomposition of the linear fundamental solution $F_t$ of (\ref{Eq: linear}) along the foliations generated by $E_1 = P \, (\R^k\times \{0\} )$ and  $E_2 = P \, (\{0\} \times \R^l)$.

\eop

Using the notation in the proof of last proposition, the decomposition of $F_t= \eta_t \circ \psi_t$ above are such that $\eta_t$ lies in the group $P\, G_H \, P^{-1}$ and $\psi_t$ lies in $P\, G_V\,  P^{-1}$.

\section{Attainability and alternate decomposition}

In this section, we introduce a technique to rescue a sort of decomposability of flows of diffeomorphisms along pair of foliations which otherwise would be topologically impossible to be performed. This is a decomposition in time intervals along the trajectories in the sense of: stop the decomposition close to the boundary of the set of diffeomorphisms where local decomposition no longer exists; and restart, from the identity, with another couple of vertical-horizontal decomposition. 
This succession of dual decomposition, vertical diffeomorphism composed with horizontal diffeomorphism, represented by $(\mathcal{H}\mathcal{V})$, will be called an \textit{alternate decomposition}. So, typically, this decomposition has the alternating structure $(\mathcal{H}\mathcal{V})\circ \cdots (\mathcal{H}\mathcal{V})$. Note that it is not relevant if the last term on the left hand side is  $\mathcal{H}$ or $\mathcal{V}$ since ending with $\mathcal{V}$ means that the omitted $\mathcal{H}$ part is the identity. The technique applies essentially in continuous flows of diffeomorphisms, or at least with jumps occurring away from sets of undecomposable diffeomorphisms, see Remark \ref{Remark: alternate for small jumps} below.

\subsection{Attainability index and topological Obstruction}

In many interesting pairs of foliations, given an initial condition $x_0\in M$, there might exist a set of points in $M$ which one cannot  reach by a vertical trajectory concatenated with a horizontal path. We investigate here the possibility of reaching these points allowing concatenation with a finite number of alternating vertical and horizontal paths, see Examples (\ref{Exemplo: aranha}) and (\ref{Exemplo: aranha infinita}) below. This topological restriction to accessibility represents also, obviously, an obstruction for the decomposition of a dynamics given by a  continuous family of diffeomorphisms $\varphi_t$ which, say, send $x_0$ into a non-accessible point. This leads us to the following concept:

\begin{definition}
\emph{
The set of $k$-\emph{attainable points} from $x \in M$ with respect to the pair of foliation $(M, \mathcal{H}, \mathcal{V})$ is the set given by the composition of saturations }
\begin{eqnarray}
\mathcal{A}^k(x) = \underbrace{\cdots \mathcal{H}(\mathcal{V}(\mathcal{H}\mathcal{V}(x))) }_{2k \mbox{ times}}.
\end{eqnarray}
\emph{
In other words, we have $k$ steps of the pair of  saturations $(\mathcal{H} \circ \mathcal{V})(\cdot)$.}
\end{definition}

\bigskip

Note that $\mathcal{A}^k(x)$ is horizontally saturated for all $k \in \mathbb{N}$ and for all $x \in M$. If a diffeomorphism  $\phi_t(x)$ is decomposable (in the sense of $(\ref{fc1}))$ in a neighbourhood of $x$, then $\phi_t(x) \in \mathcal{A}^k(x)$ for $k \in \mathbb{N}$. The converse is not true: rotations of $\pi/2$ in $\R^2$, with the Cartesian pair of foliations is a counterexamples, where $\mathcal{A}^1(x)\equiv \R^2$ for all $x\in \R^2$. Hence,  $k$-attainability is a topological obstruction to the decomposition of a diffeomorphism 
(in the sense of $(\ref{fc1}))$. 

\begin{proposition} \label{Prop: attain are open}

Given a biregular foliated space $(M, \mathcal{V}, \mathcal{H})$, the attainable sets  $\mathcal{A}^k (x)$ are open sets for all $x\in M$ and $k\in \N$.

\end{proposition}

\proof
Consider initially $k=1$ and a point $y\in \mathcal{A}^1 (x)$. By definition, there exists at least one point  $z\in \mathcal{H}(y) \cap \mathcal{V}(x) $. By uniform transversality (Theorem \ref{Thm: Uniform transversality}), there exists an open set $z\in D_1 \subset \mathcal{V}(x)=\mathcal{V}(z)$ which is sent diffeomorphically to an open set $y\in D_2\subset \mathcal{V}(y)$ along the same horizontal leaves. Using a local biregular chart at $y$ we conclude that the horizontal saturation of $D_2$ contains an open neighbourhood of $y$, hence $y$ is in the interior of $\mathcal{A}^1 (x)$.

For $k\geq 2$ one just has to write 

\[
    \mathcal{A}^k (x)= \bigcup_{y\in \mathcal{A}^{k-1} (x)} \mathcal{A}^1 (y).
\] 
The result follows by induction.

\eop

\begin{proposition}
Given a biregular foliated space $(M, \mathcal{V}, \mathcal{H})$, if $M$ is connected then $M=\cup_{k\in \N}\mathcal{A}^k(x)$ for all $x\in M$.
\end{proposition}

\proof Indeed, we only have to prove that $\cup_{k\in \N}\mathcal{A}^k$ is a closed set. Suppose that there exists a point  $x\in \partial \cup_{k\in \N}\mathcal{A}^k$. There exists a local biregular chart in a neighbourhood of $x$ which is mapped into an open rectangle in $\R^r \times \R^{n-r}$. An infinite number of points of $\cup_{k\in \N}\mathcal{A}^k$ are also mapped in this open rectangle. Trivially, these points can also reach $x$ with just one more step: vertical and horizontal trajectory (say, if $x\in \partial \mathcal{A}^k$, then $x\in \mathcal{A}^{k+1}$). We conclude that $x \in \cup_{k\in \N}\mathcal{A}^k$ hence this set is open and closed in $M$.

\eop

It is particularly interesting when one can reach the whole manifold in a finite number of steps. This leads us to the following definition:

\begin{definition}\emph{\label{def: attainability}
We define the {\it index of attainability} of the bifoliated space $(M,\mathcal{H},\mathcal{V})$ at $x \in M$ as the positive integer
\begin{eqnarray}
I_A(x, \mathcal{H},\mathcal{V}) = \mbox{min}\{k \in \mathbb{N}; \mathcal{A}^k(x) = M \},
\end{eqnarray}
when it exists. Otherwise we say that $I_A(x, \mathcal{H},\mathcal{V})=\infty$. In other words, the attainability index at $x \in M$ is the minimal number of concatenation of horizontal and vertical paths in such a way that any point on the manifold $M$ is attainable from $x$. }
\end{definition}

\begin{definition}\emph{
We define the $k$-co-attainable set of $x \in M$ with respect to $(M,\mathcal{H}, \mathcal{V})$ as the intersection
\begin{eqnarray}
\mathcal{C}^k(x) = \underbrace{\mathcal{H}\circ \mathcal{V}\circ \mathcal{H} \circ \cdots \circ \mathcal{V}(x)}_{2k \mbox{ times}} \cap  \underbrace{\mathcal{V}\circ \mathcal{H}\circ \mathcal{V} \circ \cdots \circ \mathcal{H}(x)}_{2k \mbox{ times}}.
\end{eqnarray}
Using the same argument of Proposition \ref{Prop: attain are open}, a point $y \in \mathcal{C}^k(x)$, if $y \in \mathcal{A}^k(x)$ and $x \in \mathcal{A}^k(y)$.}
\end{definition}

 
 Given $x \in M$ 
 there exists a natural sequence of open sets:
 
\begin{eqnarray}
\mathcal{A}^1(x) \varsubsetneq \mathcal{A}^2(x) \varsubsetneq \mathcal{A}^3(x) \varsubsetneq \cdots \varsubsetneq \mathcal{A}^{I_A(x,\mathcal{H},\mathcal{V})}(x)= \mathcal{A}^{I_A(x,\mathcal{H},\mathcal{V})+1}(x) = \ldots = M. 
\end{eqnarray}
Note that since the boundary  $\partial \mathcal{A}^i(x) \subset \mathcal{A}^{i+1}(x)$ then the closure $ cl \mathcal{A}^i(x) \subset \mathcal{A}^{i+1}(x)$. Moreover, since $ \mathcal{A}^i(x)$ are open, then if $M$ is compact, then ${I_A(x,\mathcal{H},\mathcal{V})}(x)$ is finite for all $x\in M$.

\begin{example}
\label{Exemplo: aranha} \emph{ 
Consider $\displaystyle M = \R^{2*} =\mathbb{R}^2\setminus\{(0,0)\}$ with the horizontal foliation given by the union of leaves:
$$
\mathcal{H} = \bigcup_{\alpha \in \R } \left\{ (x,y)\in \R^{2*}; xy = \alpha \right\}
$$
and let the vertical foliation  $\mathcal{V}$ be the rotation by $\pi/4$ of the leaves of $\mathcal{H}$. We represent $\mathcal{H}$ and $\mathcal{V}$ in Figure 1 by grey and blue curves, respectively. For any point $p$ in the diagonal $\{(x,x), x > 0\}$, we have  $\mathcal{A}^1(p) = \{ (x,y)\in \R^{2*}; y+x>0\}$, $\mathcal{A}^2(p)  = \R^2 \setminus {(x,x,), x<0}$ and $\mathcal{A}^3(p)  = M $. I.e. $I_A(p,\mathcal{H},\mathcal{V}) = 3$, and
\begin{eqnarray*}
\mathcal{C}^1(p) &=& \{(x,y)\in \R^{2*};\  x>0, y>0\}\\
\mathcal{C}^2(p) &=& \{(x,y)\in \R^{2*};\  y<0 \mbox{ or } x<0 \} \\
\mathcal{C}^3(p) &=& M.
\end{eqnarray*} }
\end{example}

\begin{figure}[ht]
    \centering
\includegraphics[scale=0.50]{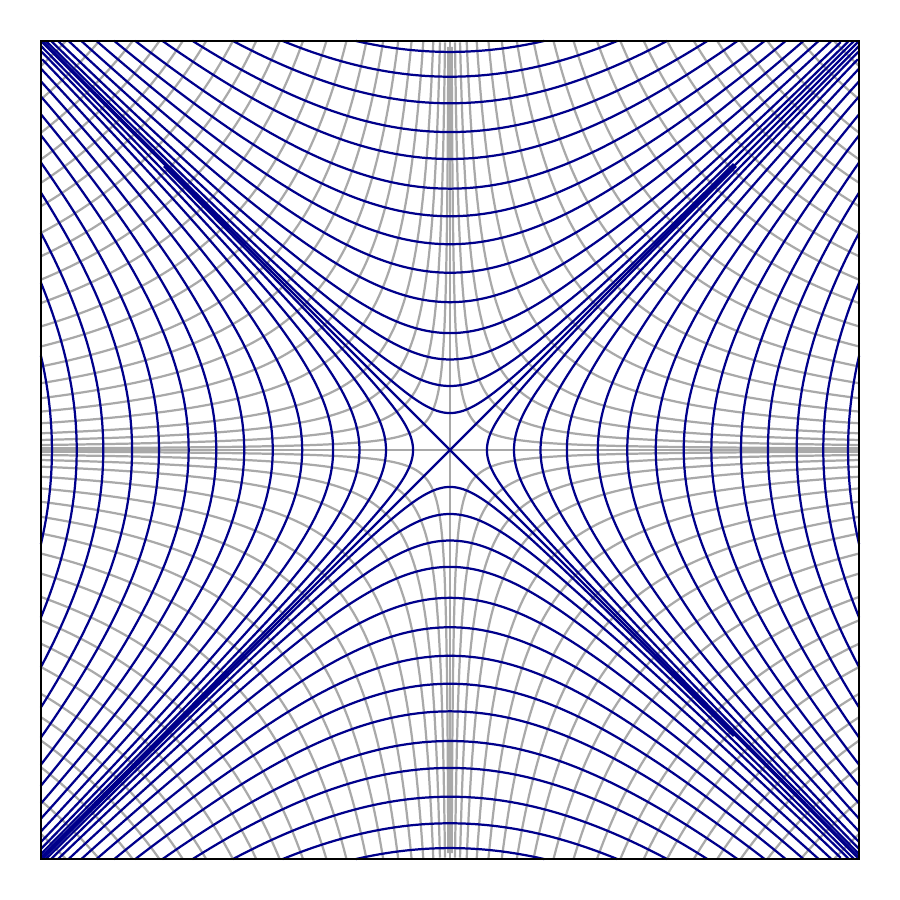}
    \caption{Example of a bifoliation with $I_A(p,\mathcal{H}, \mathcal{V})=3$ } \label{figura1}
\end{figure}

\eop

In many interesting cases the index of attainability is infinite:

\begin{example} \emph{For $M = \mathbb{R}^2 $ consider the horizontal foliation given by the disjoint union of leaves
$$
\mathcal{H} = \bigcup_{c\in \R} \big\{\{(x,y) \in \mathbb{R}^2; y =  -\Big|\sec x \Big| + c, \mbox{ if } x\neq r\pi \Big\}  \ \    \dot{\bigcup} \ \ \Big\{x = r\pi; \ r \in \mathbb{Z}\Big\},
$$
and vertical foliation given by 
$$
\mathcal{V} =  \bigcup_{c\in \R} \big\{\{(x,y) \in \mathbb{R}^2; y =  -\Big| \csc x \Big| + c, \mbox{ if } x\neq r\pi +\frac{\pi}{2}\Big\}  \ \    \dot{\bigcup} \ \ \Big\{x = r\pi + \frac{\pi}{2}; \ r \in \mathbb{Z}\Big\},
$$
see Figure \ref{figura2}. In this case, note that for $p = (0,1)$:  
\begin{eqnarray*}
\mathcal{A}^1(p) &=& \Big\{(x,y) \in \mathbb{R}^2; -\frac{\pi}{2} \leq x \leq \frac{3\pi}{2} \Big\} \\
\mathcal{A}^2(p) &=& \Big\{(x,y) \in \mathbb{R}^2; -\frac{3\pi}{2} \leq x \leq \frac{5\pi}{2} \Big\}.  
\end{eqnarray*}
In general,
\begin{eqnarray*}
\mathcal{A}^n(p) &=& \Big\{(x,y) \in \mathbb{R}^2 ; -n\pi + \frac{\pi}{2} \leq x \leq n\pi +\frac{\pi}{2} \Big\}.
\end{eqnarray*}
Hence  
\[
M = \bigcup_{n=1}^{\infty}\mathcal{A}^n(p),
\]
and the index of attainability for any initial point $p$ is given by  $I_A(p,\mathcal{H}, \mathcal{V}) = \infty$. }

\eop

\label{Exemplo: aranha infinita}

\begin{figure}[ht]
    \centering
    \includegraphics[scale=0.43]{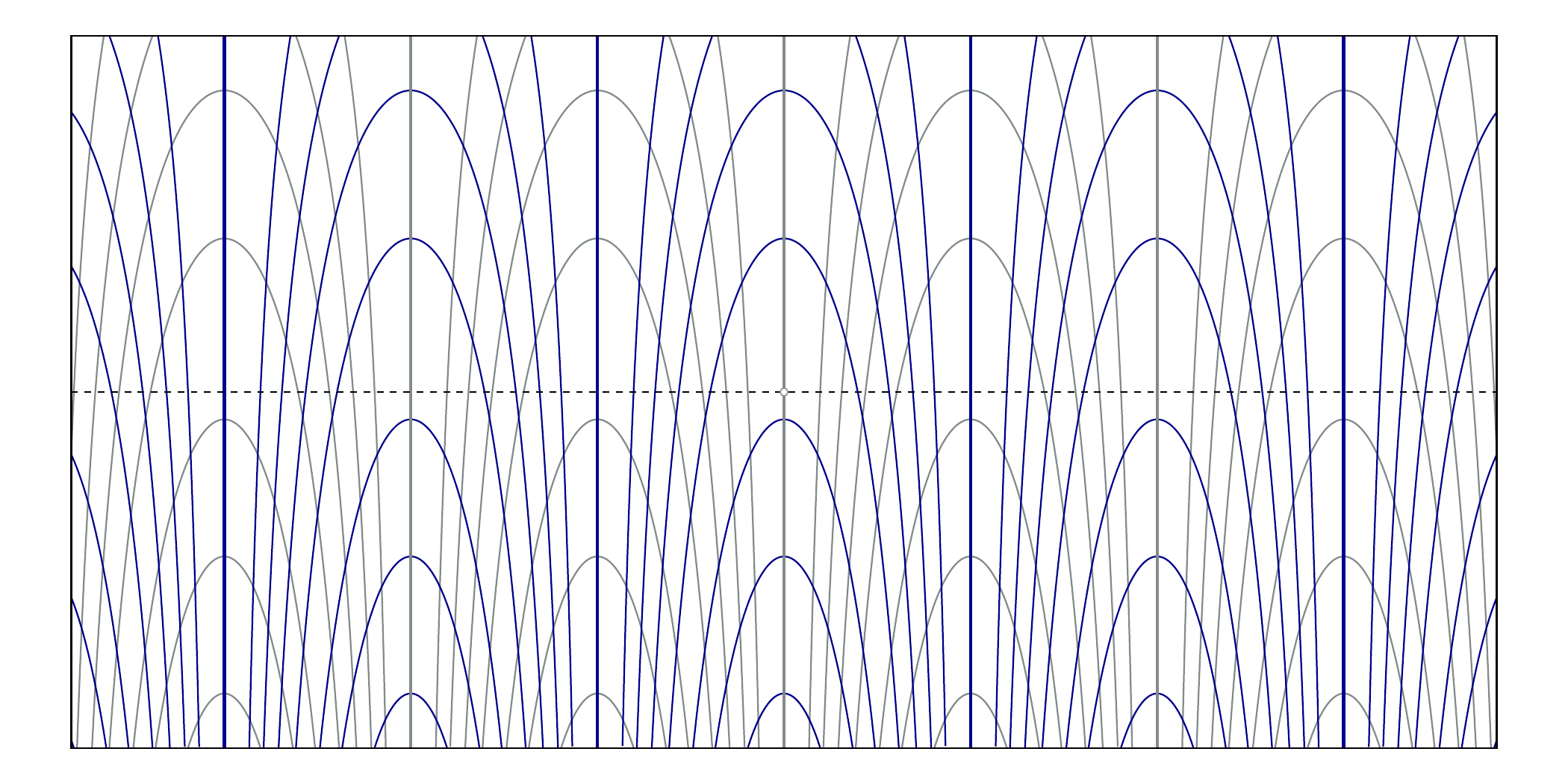}
    \caption{Example of a  bifoliation with $I_A(p,\mathcal{H}, \mathcal{V})= \infty$ for any point $p$. }
    \label{figura2}
\end{figure}

\end{example}

\begin{proposition} \label{Prop: Attan}
    Let $M$ be a connected 
    manifold. Consider  $p \in M$ and $k \in \mathbb{N}$ such that horizontal and vertical saturations commute in the sense that $\mathcal{A}^k(p) = \underbrace{\mathcal{V}\circ \mathcal{H} \circ  \ldots  \mathcal{V}  \circ \mathcal{H}(p)}_{2k \mbox{ times}}$, then $\mathcal{A}^k(p) = M$, i.e. the index of attainability $I_A(p,\mathcal{H}, \mathcal{V})= k$. In particular, in this case, $M= \mathcal{C}^{I_A(p,\mathcal{H},\mathcal{V})}(p)
$.

\end{proposition}

\proof Since  $ \mathcal{A}^k(p)$ is an open set and $M$ is connected, we only have to prove that $ \mathcal{A}^k(p)$ is closed. By definition $ \mathcal{A}^k(p)$ is horizontally saturated, hence its boundary $\partial\mathcal{A}^k$ is also horizontally saturated by uniform transversality (Thm.  \ref{Thm: Uniform transversality}). Hence, given a point  $x \in \partial\mathcal{A}^k$, since  $\mathcal{H}$ and $\mathcal{V}$ are complementary, there exists a point $y$ in the intersection $\mathcal{V}(x) \cap \mathcal{A}^k (p)$. The hypothesis says that $\mathcal{A}^k(p)$ is also vertically saturated, hence, the point in the boundary  $x\in \mathcal{V}(y)$ must also be in  $\mathcal{A}^k (p)$. The last statement follows straightforward by the definition. 

\eop

\subsection{Alternate decomposition}

In general, let $(\phi_t)_{t}$, be a continuous family of global diffeomorphisms on $M$, with $\phi_0 = Id$ and $0 \leq t\leq T$. Fix an initial condition $p \in M$. The single decomposition  $\phi_t= \eta_t\circ \psi_t$ where $\eta_t$ is horizontal and $\psi_t$ is vertical for all time $t \in [0,T]$ in general does not exist. In fact, we have the following topological characterization of decomposability of a flow with just one pair horizontal-vertical:

\begin{theorem}
\label{theorem: Melo}
Suppose that $(M, \mathcal{H}, \mathcal{V} )$ is transversely orientable for the horizontal foliation. Then a family of diffeomorphisms $\phi_t$ is globally decomposable (in the sense of Remark \ref{characterization of decomposition}) for all $0 \leq t \leq T$, if and only if, it preserves  transverse orientation.
\end{theorem}

\proof 
For a proof, see  \cite[Thm. 2.5]{Melo1}.

\eop

 This restriction motivates the alternate decomposition of equation (\ref{fc1}), whose existence is treated in the sequel. Before that, next Proposition, give us a simple necessary condition of its existence.

\begin{proposition}
Suppose $\phi_t$ is a flow of diffeomorphism which can be decomposed as 
\begin{eqnarray*}
\phi_t(x) = \big(\eta^k_{t\vee s_{k-1}}\circ\psi^k_{t\vee s_{k-1}} \big) \circ \ldots \circ \big( \eta^2_{s_2}\circ\psi^2_{s_2}\big) \circ \big(\eta^1_{s_1}\circ\psi^1_{s_1}\big),
\end{eqnarray*}
up to a time $\tau>0$, where $\psi^i$ and $\eta^i$ are purely vertical and horizontal diffeomorphisms, respectively, and $0< s_1 < s_2< \ldots <s_k < s_{k+1}<  \ldots  < \tau $ is a non-decreasing sequence of times.  Then, for $0\leq t<\tau$ and $i=1,\ldots, k$ we have that   $\psi^i_t(x)\in \mathcal{V}(x)\cap \mathcal{H}\left( \eta_t^i\circ \psi_t^i(x)\right)$ and $\eta^i_t(y) \in \mathcal{H}(y)\cap \eta_t^i\circ \psi_t^i\left(\mathcal{V}(y)\right)$.
\end{proposition}

\proof
In fact, for the first statement, note that $\psi_t^i(x)\in \mathcal{V}(x)$, for $i=1,\ldots, k$ and $x$ in the appropriate domain. In addition, $\eta_t$ preserves horizontal leaves for all $t<\tau$, then since $\mathcal{H}\left(\xi^i_t\circ\psi^i_t(x)\right) = \mathcal{H}\left(\psi^i_t(x)\right)$, it implies that $\psi^i_t(x)\in\mathcal{H}\left(\xi_t^i\circ\psi^i_t(x)\right)$. For the second statement, observe that $\eta_t^i(y) \in \mathcal{H}(y)$ and $y\in\mathcal{V}(y)$, therefore $\eta_t^i(y)\in \eta_t^i\circ \psi_t^i\left(\mathcal{V}(y)\right)$.

\eop

\bigskip

\begin{theorem}
\label{teoremao}
Let $(\phi_t)_{t \in [0,T]}$  be a family of diffeomorphisms on $M$ depending continuously on $t$, with $\phi_0 = Id$ and $T\in \R \cup \{+\infty\}$. For any fixed $x \in M$,  there exists a non-decreasing sequence of times $0= s_0< s_1 < s_2< \ldots <s_r =s_{r+1}= \ldots  = T $, $(r\in \N \cup \{\infty\})$ such that for all $t \in [0,T]$, there exists a neighbourhood  $U_x$ of $x$, where we have the following foliated decomposition:
\begin{eqnarray*}
\phi_t(x) = \big(\eta^k_{t}\circ\psi^k_{t} \big) \circ \ldots \circ \big( \eta^2_{s_2}\circ\psi^2_{s_2}\big) \circ \big(\eta^1_{s_1}\circ\psi^1_{s_1}\big),
\end{eqnarray*}
for  $t\in [s_{k-1}, s_k]$, with $\eta^k_{ s_{k-1}}= Id$ and $ \psi^k_{ s_{k-1}}= Id$. Here $\eta^j$ and $\psi^j$ are horizontal and vertical diffeomorphisms respectively for all $j \in 1, \ldots, k$.
\end{theorem}

\proof
For each $x_0 \in M$, we consider a neighbourhood $U_{x_0}$ of $x_0$ and a coordinate system \\ $\displaystyle{\Sigma_{x_0}: U_{x_0} \subset M \longrightarrow \mathbb{R}^k\times \mathbb{R}^{n-k} }$; with respect to it, we write $\phi_s(x) = (\phi^1_s(x,y), \phi^2_s(x,y))$, for $x \in U_{x_0}$, $s>0$ and we define 
\begin{eqnarray*}
s_1 = \inf \left\{s \in [0,T]; \det \frac{\partial \phi^2_t}{\partial y} = 0 \right\} - \epsilon_1,
\end{eqnarray*}
for $\epsilon_1 > 0$ small enough. For $s<s_1$, $\phi_s(x_0)$ preserves transverse orientation. Applying a local version of Theorem (\ref{theorem: Melo}), it follows that in $U_{x_0}$, $\phi_s(x)$ has a foliated decomposition 
\begin{eqnarray*}
\phi_s(x) = \eta^1_s \circ \psi^1_s(x),
\end{eqnarray*}
where $\eta^1_s$ and $\psi^1_s$ are horizontal and vertical continous family of diffeomorphisms respectively. If $s_1 = T$, the proof stops, otherwise, for $s>s_1$, we can write $\phi_s(x)$ as 
\begin{eqnarray*}
\phi_s(x) &=& \phi_{s_1,s}\circ \phi_{s_1}(x) \\
&=& \phi_{s_1,s} \circ (\eta^1_{s_1} \circ \psi^1_{s_1})(x),
\end{eqnarray*}
where $\phi_{s_1,s}=\phi_{s}\circ \phi_{s_1}^{-1}$.
Taking $\phi_{s_1,s}$ sufficiently close to the identity, the decomposition can be performed again. More precisely, set 
\begin{eqnarray*}
s_2 = \inf \left \{s \in [s_1, T]; \det \frac{\partial \phi^2_{s_1,s}}{\partial y}(\phi_{s_1}(x_0)) = 0 \right \} - \epsilon_2,
\end{eqnarray*}
for a sufficiently small $\epsilon_2 > 0$. For $s_1 < s < s_2$, we apply Theorem (\ref{theorem: Melo}) again and rewrite $\phi_s(x)$ as 
\begin{eqnarray*}
\phi_s(x) = (\eta^2_{s} \circ \psi^2_{s}) \circ (\eta^1_{s_1} \circ \psi^1_{s_1})(x),
\end{eqnarray*}
where $\eta^2$ and $\psi^2$ are horizontal and vertical diffeomorphisms respectively and the domain is restricted if necessary for the second decomposition. If $s_2 =T$, the proof stops, otherwise we repeat the argument recursively as many times as necessary until $t\geq s_k$. For a fixed $t$, parameter $k$ does not go to infinity since in the compact interval $[0,T]$ the derivatives of $\phi_t$ are bounded and by compactness, the trajectory cross only a finite number of boundaries of attainable sets $\mathcal{A}^i(x_0)$.


\eop

\begin{proposition}
\label{propdet}
Suppose that $(M,\mathcal{H}, \mathcal{V})$ is transversely orientable for the horizontal foliation. If $\phi_s(p)$ approaches the boundary $\partial \mathcal{A}^k(p)$ of the $k$-attainable set, the determinant of  \  $\frac{\partial \phi_s^2}{\partial y}(\phi_{s_{k-1}}(p))$ goes to zero.
\end{proposition}

\proof In fact, by definition of $k$-attainable sets, their boundaries indicate the topological limit of the attainability with vertical trajectories concatenated with horizontal trajectories. Hence it is also an obstruction for the decomposition of a flow.  Cf. Equation (\ref{prop1}), decomposition exists if and only if the derivative $\det \frac{\partial \phi^2_{s_0}(x,y)}{\partial y}\neq 0$. 

\eop

\bigskip


The next example shows that Theorem (\ref{teoremao}) may hold even if ($M,\mathcal{H},\mathcal{V}$) is not transversely orientable for the horizontal foliation.

\bigskip

\begin{example}
\emph{
Let $M = [0,1]^3/\sim $, where $\sim$ is the identification of the following faces of the cube $[0,1]^3$: }
\begin{eqnarray}
(x,0,z) \sim  (1-x,1,1-z),
\end{eqnarray}
\emph{such that the section $\left(x,y, \frac{1}{2}\right) \cap [0,1]^3$ turns into a M\"obius strip $S$. Note that $M$ is a tubular neighbourhood of $S$. In this context, the horizontal and vertical  foliations $\mathcal{H}$ and $\mathcal{V}$ are given by the image of the horizontal and vertical plaques respectively.  It is worth mentioning that ($M\setminus S, \mathcal{H}$) is transversely orientable, but $(M,\mathcal{H})$ is not.
Consider a complete family of diffeomorphisms given by $\phi_t(x,y,z) = (x,y+t,z)$. In this case, $\phi_t$ is a horizontal flow with respect to the pair of foliation ($\mathcal{H}, \mathcal{V}$), hence it can be decomposed as $\phi_t = \eta_t \circ \psi_t$, where $\eta_t = \phi_t$ and $\psi_t = Id$, for small $t>0$. Around the non-transversely orientable leaf $S$, $\phi_t$ has an alternate decomposition  $\displaystyle{\phi_t(x_0) = \big(\eta^k_t\circ\psi^k_t \big) \circ \ldots \circ \big( \eta^2_{s_2}\circ\psi^2_{s_2}\big) \circ \big(\eta^1_{s_1}\circ\psi^1_{s_1}\big)(x_0)}$, where $\displaystyle{s_i \in \{(2k+1)2\pi, \hbox{with} \  k \in \mathbb{Z}\}}$, for a local biregular coordinate system in a neighbourhood of an initial condition $x_0 \in S$. Each pair $\eta^j_t \circ \psi^j_t$ is given by the projection of the two reverting orientation diffeomorphisms  $\eta^j_t(x,y,z) = (y,x,z)$ and $\psi^j_t(x,y,z) = (x,y,1-z)$, since $\phi_t$ reverses the orientation of both vertical and horizontal components just before $t = s_i$.  In the manifold ($M \setminus S, \mathcal{H}$), the decomposition is guaranteed by Theorem (\ref{theorem: Melo}).}
\end{example}

 \eop

We end this Section with a particular case of  alternate decomposition for classical stochastic flows of diffeomorphisms. 

\begin{theorem}
\label{teorem do caso estocastico}
Let $\varphi_t$ be the stochastic flow of (local) diffeomorphisms generated by the SDE (\ref{eq6}). There exists a non-decreasing sequence of stopping times $0 = t_0 < t_1 < t_2 < \ldots < t_r = t_{r+1} = \ldots = a$ such that, locally, $\varphi_t$ is alternately decomposable as 
\begin{eqnarray}
\varphi_t(\omega,x) = \big(\eta^k_{t_{k-1},t}\circ\psi^k_{t_{k-1},t} \big) \circ \ldots \circ \big( \eta^2_{t_1,t_2}\circ\psi^2_{t_1,t_2}\big) \circ \big(\eta^1_{t_1}\circ\psi^1_{t_1}\big)(\omega, x),
\end{eqnarray}
where $\eta^j$ and $\psi^j$ are horizontal and vertical diffeomorphisms respectively for all $j \in \mathbb{N}$, $t \in [t_{k-1}, t_k]$ and $\omega \in \Omega$.
\end{theorem}

\proof Theorem (\ref{teoremao}) guarantees that the alternate decomposition holds pathwise. Hence, we only have to  consider a sequence of stopping times $(t_i)_{i \in \mathbb{N}}$ which satisfies (recursively): 
\begin{eqnarray}
t_i(\omega,x) = \inf \left\{t \in [t_{i-1}(\omega,x),a); \det \frac{\partial \varphi^2_{t_{i-1},t}}{\partial y}(u_i) = 0 \right\} - \epsilon_i,
\end{eqnarray}
where $u_i = \varphi_{t_i}(x)$, for $\epsilon_i>0$ small enough, $\omega \in \Omega$. By the cocycle property for stochastic flows, we can write $\varphi_t$ initially as a composition of diffeomorphisms which are all sufficiently close to the identity; after that, we decompose each of these components according to the foliations, i.e.
\begin{eqnarray*}
\varphi_t(\omega,x) &=& \varphi_{t_k,t}(\theta_{t_k}(\omega),u_k) \circ \ldots \circ \varphi_{t_1,t_2}(\theta_{t_1}(\omega),u_1) \circ \varphi_{t_1}(\omega,x) \\
&=& \big(\eta^k_{t_{k-1},t}\circ\psi^k_{t_{k-1},t} \big) \circ \ldots \circ \big( \eta^2_{t_1,t_2}\circ\psi^2_{t_1,t_2}\big) \circ \big(\eta^1_{t_1}\circ\psi^1_{t_1}\big)(\omega, x),
\end{eqnarray*}
where $\theta_t$ is the canonical shift operator on the probability space and $\eta^j$, $\psi^j$ are horizontal and vertical stochastic flows of diffeomorphisms.

\eop

\begin{remark} \label{Remark: alternate for small jumps} \emph{Although along this section we have treated  alternate decomposition only on continuous flows of diffeomorphisms, it can also be applyed to jumping systems. In fact, if the jumps are only inside $\mathcal{A}^{1}$, or more generally, if the trajectories cross boundaries  of attainable sets  $\mathcal{A}^{i}$ only at times of continuity with respect to $t$, the alternate decomposition can also be performed.  }
\end{remark}

\section{Principal fibre bundles over homogeneous spaces}


Consider a connected Lie group $G$ with a closed subgroup $H$ and let $\mathfrak{g}$ and $\mathfrak{h}$ represent their corresponding Lie algebras of right invariant vector fields. The action of $G$ on $H$ is given by left translation $gH$,  for any $g \in G$ and the orbits generate the homogeneous space $M:=G/H$, as described in the literature e.g. \cite{kobayashi}. There exists a principal fibre bundle characterized by the canonical projection $\pi: G \rightarrow M$.  For any element $A\in \mathfrak{g}$ consider the right invariant SDEJ:
\begin{eqnarray}
\label{equation: principal fibre bundle}
d\, g_t =  A  g_t\diamond dZ_t.
\end{eqnarray}

Consider a connection $\omega$ in the principal fibre bundle $\pi: G \rightarrow M$. We construct the decomposition of flow according to the vertical subspace (involutive) and the horizontal subspace established by this connection. The solution flow (global in $G$ up to lifetime of $Z_t$) is given by the left action:
\[
\varphi_t (x)= g_t x,
\]
where $g_t= \exp\{ A\, Z_t \}$. The distributions $\Delta_{\mathcal{H}}$ and $\Delta_{\mathcal{V}}$ in the tangent space $TG$ are defined by the horizontal subspace with respect to the connection $\omega$ and the tangent to the fibres $gH$ (involutive). In order to decompose the flow $\varphi_t$ as outlined in Theorem (\ref{disc version}),  it is necessary to determine the vector fields  $V$ and $h$ as in equation (\ref{horizontal vector fields}) in the proof of the Theorem, i.e.:
\[
Ax := h + V(\eta, x).
\]
Elements $\eta \in \mathrm{Diff}(\Delta_{\mathcal{H}}, G)$ can be expressed pointwise (with respect to $x\in G$) as a left action by elements of $G$ at $x$. This is a vertical component preserving action, i.e. $g_* \Delta_{\mathcal{V}} = \Delta_{\mathcal{V}}$ for all $g\in G$. Hence: 
\[
V(x)= \omega (Ax)^*  \ \ \ \     \mbox{ and }  \ \ \ \ h=Ag- \omega (Ax)^*.
\]
From (\ref{jump horizontal}) and (\ref{jump vertical}) it follows that each component of the decomposition $\varphi_t(\cdot)= \eta_t \circ \psi_t (\cdot) $ are given by:
\begin{equation} \label{Eq: eta para fibrado}
d\, \eta_t = R_{\eta_t*} (A\eta_t(\cdot) - \omega (A\eta_t (\cdot))^*)
\end{equation}
and 
\begin{equation} \label{Eq: psi para fibrado}
d\, \psi_t = \mathrm{Ad}(\eta_t) \ \omega (A \eta_t (\cdot))^*.
\end{equation}

Let denote by $g^{\mathcal{H},x}_t \in G$ the general semimartingale in $G$ such that $g^{\mathcal{H},x}_t x $  is the horizontal lift of $\pi(g_t\, x)$  at the point $x$, i.e. $g^{\mathcal{H},x}_t x $ is a horizontal càdlàc path and $g^{\mathcal{H},x}_t x = g_t\, x\, v_t$ for some $v_t \in H$. Using this notation and fixing the action at a point $x\in G$, the previous equations simplify to well known finite dimensional equations (in $G$). This concept is detailed in the next subsection.

\begin{proposition}
\label{Prop: decomposition in fibre bundles}
Consider the decomposition of the solution flow  $\varphi_t (\cdot)=  \eta_t \circ \psi_t (\cdot)$  for equation (\ref{equation: principal fibre bundle}) in accordance with the horizontal and vertical distribution of the fibre bundle as described in Theorem (\ref{disc version}). Then, at each point $x\in G$, the first component can be written as the left action: 
\[
\eta_t(x)= g^{\mathcal{H},x}_t x, 
\]
and the second component can be written as the right action:
\[
\psi_t(x)= x \, h_t
\]
where $h_t=
x^{-1}\,(g^{\mathcal{H},x}_t)^{-1} \, g_t \, x.
$
\end{proposition}
\proof The proof of the first equation follows straightforward when equation (\ref{Eq: eta para fibrado}) is applied  at a fixed initial condition $x\in G$. Concerning the second equation mentioned in the statement, one sees that $(x^{-1}\,(g^{\mathcal{H},x}_t)^{-1} \, g_t \, x )\in H$ by definition of the horizontal lift: $g^{\mathcal{H},x}_t x = g_t\, x\, v_t$ for some $v_t \in H $. One checks that it solves (\ref{Eq: psi para fibrado}) at a fixed point $x$.

\eop

\bigskip

\subsection{Trivial fibre bundles} 

As a particular case, Let  $\pi : G\times H \rightarrow H$ be a trivial principal fibre bundle with structural group $H$, with connected Lie groups $G$ and $H$. The trivial connection here is given by $\omega_{(x,y)}(g_t', h_t')= y^{-1} h_t'\in \mathfrak{h}$. Consider a right invariant SDEJ in $G \times H$:
\[
d (x_t, y_t) = (A\times B) \  (x_t, y_t) \diamond dz_t.
\]
Where $A\in \mathfrak{g}$ and $B \in \mathfrak{h}$, the Lie algebras of $G$ and $H$ respectively, with an initial condition $(x_0, y_0)$. By the fact that, in this case, the connection  is invariant by left action of $G\times \{1d\}$, the factor $g^{\mathcal{H},x}_t\in G \times H$ of Proposition 
(\ref{Prop: decomposition in fibre bundles}) does not depend on $(x,y)$. The trivial components of the decomposition are recovered. Indeed, a global decomposition is achieved where the first component results from the left action:
\[
\eta_t(\cdot, \cdot)= (\exp (Az_t), 1d) (\cdot,\cdot). 
\]
And the second (vertical) component is given in terms of the right action:
\[
\psi_t(\cdot ,\cdot)= (\cdot , \cdot ) (1d, h_t),
\]
where $h_t=
y^{-1}\, \exp (Bz_t) \, y 
$.





\subsection{Jump dynamics on reductive homogeneous spaces}

 \noindent A homogeneous space $M=G/H$ is said to be reductive if the Lie algebra $\mathfrak{g}$ contains a subspace $\mathfrak{n}$, such that $\mathrm{Ad}(H)(\mathfrak{n}) \subset \mathfrak{n}$ and $\mathfrak{g}$ can be written as the direct sum $\mathfrak{g} = \mathfrak{h} \oplus \mathfrak{n}$. It is worth mentioning that each fibre $\pi^{-1}(x)$ is diffeomorphic to $H$. A similar decomposition was considered by Li  \cite{Li} in the context of standard Brownian motion.  

 \begin{corollary}
 \label{corollary: ITô-Ventzel in Lie groups}
 \emph{
 For $S \in C_b^2(G, \mathcal{L}(\mathbb{R}^d, \mathfrak{h}))$ and $Y \in C_b^2(G, \mathcal{L}(\mathbb{R}^d, \mathfrak{n}))$. Let $Z_t$ be a general semimartingale and $\psi_t$ and $\eta_t$ be solutions associated with the Marcus differential equations $d\psi_t = S^*(\psi_t) \diamond dZ_t$ and $d\eta_t = Y^*(\eta) \diamond dZ_t$. Then,} 
 \begin{eqnarray}
 d(\eta_t \psi_t) = R_{\psi_t^*}d\eta_t + (L_{\psi_t^{-1}}d\psi_t)^*(\eta_t \psi_t).
\end{eqnarray}
 \end{corollary}
 \proof
 By Corollary \ref{corol1}, it follows that 
 \begin{eqnarray*}
 d(\eta_t \psi_t)  &=& R_{\psi_t*}Y^*(\eta_t)dZ_t + L_{(\eta_t \psi_t)^*}L_{\psi^{-1}*}S^*(\eta_t \psi_t)dZ_t \\
 &=& R_{\psi_t*}d\eta_t + (L_{\psi^{-1}}d\psi)^*(\eta_t \psi_t).
 \end{eqnarray*}

 \eop
 
 \bigskip
 
 Let $\varphi_t, t < T$, be the flow of diffeomorphism of the following canonical Marcus stochastic differential equation,
\begin{equation}
\label{Jump main eq in Lie groups}
d\varphi_t =  W^*(\varphi_t)\diamond dZ_t.
\end{equation}
Where $W$ is a element of the Lie algebra $\mathfrak{g}$.

\bigskip

 In the next theorem, we find explicit Marcus differential equations for the vertical and horizontal components of the solution $\varphi_t$. Let $\omega$ be the canonical $1$-form connection on the principal bundle $(P,G,H,\pi, M)$. By definition, $\omega(X) = 0$ for all vector field $X \in \mathfrak{n}$, and $\omega(A^*)=A$ if $A \in \mathfrak{h}$.  We suppose that the Lie algebra $\mathfrak{g}$ is reductive, therefore it can be written as the direct sum $\mathfrak{g} = \mathfrak{h} \oplus \mathfrak{n}$. Thus, the vector field $W$ can be decomposed into $W^*(g) = h^*(g)+V^*(g)$, where $h^*(g) \in \mathfrak{h}$ and $V^*(g) \in \mathfrak{n}$.

\begin{theorem}
\emph{
The solution flow $\varphi_t$ can be decomposed into $\varphi_t = \eta_t \psi_t$, such that the components $\eta_t$ and $\psi_t$ satisfies the following system of Marcus differential equations:}
\begin{eqnarray}
d\psi_t &=& V^*(\psi_t)\diamond dZ_t^i, \label{decomp1}\\
d\eta_t &=&  \left(\mbox{Ad}(\psi_t)h\right)^*(\eta_t)\diamond dZ_t. \label{decomp2}
\end{eqnarray}
\end{theorem} 
\proof
Note that $\varphi_t$ is a diffeomorphism for all $t>0$, moreover, it sends each fibre in another fibre, hence $\varphi_t$ is decomposable for all $t< T$, see e.g. \cite[Corollary 2]{Melo2}, thus the solution flow $\varphi_t$ can be rewritten as $\varphi_t=\eta_t \circ \psi_t$, where $\eta_t$ and $\psi_t$ are horizontal and vertical semimartingales respectively. By Corollary \ref{corollary: ITô-Ventzel in Lie groups}, we have that:
\begin{eqnarray}
d\varphi_t = R_{\psi_t} d\eta_t + \left(L_{\psi_t^{-1}}d\psi_t\right)^*(\varphi_t). \label{ivkdec1}
\end{eqnarray}

\noindent Applying the 1-form $\omega$ at $d\varphi_t$:
\begin{eqnarray*}
\omega(d\varphi_t) = \omega\left(L_{\psi_t^{-1}}d\psi_t\right)^*(\varphi_t) = L_{\psi_t^{-1}}d\psi_t.
\end{eqnarray*}
Hence,
\begin{eqnarray*}
L_{\psi_t^{-1}}d\psi_t &=& \omega\left( W^*(\varphi_t)\diamond dZ_t^i\right) \\
&=& \omega\left( h^*(\varphi_t)\diamond dZ_t +  V^*(\varphi_t)\diamond dZ_t\right) \\
&=&  V\diamond dZ_t. \\
\end{eqnarray*}
Therefore,
\begin{eqnarray*}
d\psi_t &=&  V^*(\psi_t)\diamond dZ_t. 
\end{eqnarray*}

\noindent Using the identity $\psi_t\circ \psi_t^{-1} = 1$ and Corollary \ref{corollary: ITô-Ventzel in Lie groups}, it follows that:
\begin{equation*}
d\psi_t^{-1} = - R_{\psi_t^{-1}} V \diamond dZ_t.
\end{equation*}

\noindent Applying Corollary \ref{corollary: ITô-Ventzel in Lie groups} once again, for $\eta_t=\varphi_t\circ\psi_t^{-1}$, we have:
\begin{eqnarray*}
d\eta_t &=& R_{\psi_t^{-1}}d\varphi_t + L_{\varphi_t}d\psi_t^{-1} \\
&=& R_{\psi_t^{-1}} \left( W^*(\varphi_t)\diamond dZ_t\right) - L_{\varphi_t}\left( R_{\psi_t^{-1}} V \diamond dZ_t\right) \\
&=&  \left[R_{\psi_t^{-1}} L_{\varphi_t} W - R_{\psi_t^{-1}}L_{\varphi_t} V\right] \diamond dZ_t \\
&=&  \left[R_{\psi_t^{-1}} L_{\eta_t}L_{\psi_t}h\right] \diamond dZ_t \\
&=&  L_{\eta_t} \mbox{Ad}(\psi_t) h \diamond dZ_t =  \left( \mbox{Ad}(\psi_t) h\right)^* (\eta_t) \diamond dZ_t.  \hspace*{10cm}
\end{eqnarray*}

\eop

\bigskip

\noindent Let $\nu_t = \pi(\varphi_t)$, we want to compute a Marcus differential equation for $\nu_t$. Using the fact that $d\pi(V^*(\varphi)) = 0$, it follows that:
\begin{eqnarray*}
d\nu_t = d\pi(d\varphi_t) &=& d\pi\left(W(\varphi_t)\right)\diamond dZ_t \\
&=&  d\pi\left(h^*(\varphi_t) + V^*(\varphi_t)\right)\diamond dZ_t \\
&=&  d\pi h^*(\varphi_t) \diamond dZ_t. \hspace*{12cm}
\end{eqnarray*}

\noindent Therefore
\begin{equation}
d\nu_t =  \bar{L}_{\eta_t*} \bar{L}_{\psi_t*} d\pi(h) \diamond dZ_t,
\label{decomp3}
\end{equation}
where $\bar{L}_a$ is the left translation on the base space $M$, for $a\in G$. Then, $\pi\circ L_a = \bar{L}_a \circ \pi$.

\begin{proposition}
\emph{
The process $\eta_t, t< T$ satisfies the equation (\ref{decomp2}), if, and only if, it is a horizontal lift of $\nu_t$.
}
\end{proposition}
\proof
Suppose that $\eta_t$ is a solution flow of equation (\ref{decomp2}). 
Since $\omega_{\eta_t}(d\eta_t)=0$, and $\psi_0\in H$, taking $x_t=\pi(\eta_t)$, it holds that:
\begin{eqnarray*}
d x_t &=&  d\pi\left( \mbox{Ad}(\psi_t^{-1})h \right)^*(\eta_t) \diamond dZ_t \\
&=&  d\pi \left(R_{\psi_t^{-1}}L_{\psi_t} h\right)^*(\eta_t) \diamond dZ_t \\
&=&  \bar{L}_{\eta_t} \bar{L}_{\psi_t} d\pi h \diamond dZ_t, \hspace*{12cm}
\end{eqnarray*}
therefore, $x_t$ satisfies equation (\ref{decomp3}), by uniqueness of solution of Marcus differential equations, we conclude that $x_t = \pi(\eta_t)$.

On the other hand, suppose that $\eta_t$ is a horizontal lift of $\nu_t$ up to a stopping time $T$. Since $\varphi_t$ is a solution of equation (\ref{Jump main eq in Lie groups}) and $\pi(\varphi_t) = \pi(\eta_t)$, then $\varphi_t$ and $\eta_t$ belong to the same fibre for all $t < T$. Therefore, there exists $C_t \in G$, such that $\eta_t C_t = \varphi_t$, for $t < T$. By Corollary \ref{corollary: ITô-Ventzel in Lie groups}: 
\[
d\eta_t = R_{C_t^{-1}}d\varphi_t + \left(L_{C_t}dC_t^{-1}\right)^*(\eta_t).
\]
We rewrite the above expression by:
\begin{eqnarray}
\label{expression: horizontal lift}
\nonumber d\eta_t &=& R_{C_t^{-1}}\left(h^*(\varphi_t)+ V^*(\varphi_t)\right) \diamond dZ_t + \left(L_{C_t}dC_t^{-1}\right)^*(\eta_t) \\
 &=& R_{C_t^{-1}}h^*(\varphi_t) \diamond dZ_t + R_{C_t^{-1}}V^*(\varphi_t) \diamond dZ_t + \left(L_{C_t}dC_t^{-1}\right)^*(\eta_t) \\
\nonumber &=&  \left(\mbox{Ad}(C_t)h\right)^*(\eta_t) \diamond dZ_t + \left(\mbox{Ad}(C_t)V\right)^*(\eta_t) \diamond dZ_t + \left(L_{C_t}dC_t^{-1}\right)^*(\eta_t).
\end{eqnarray}
Now, applying the connection $1$-form $\omega$ to the expression (\ref{expression: horizontal lift}): 
\begin{eqnarray}
\label{expression2: horizontal lift}
\nonumber 0 &=& \omega_{\eta_t}\left(R_{C_t^{-1}}V^*(\varphi_t)\right) \diamond dZ_t + L_{C_t}dC_t^{-1} \\
&=& \omega_{\eta_t}\left((\mbox{Ad}(C_t)(V))^*(\eta_t)\right) \diamond dZ_t + L_{C_t}dC_t^{-1} \\
\nonumber &=& \mbox{Ad}(C_t)V \diamond dZ_t + L_{C_t}dC_t^{-1}.
\end{eqnarray}
Here we used the fact that $\omega(d\eta_t) = 0$ and $\mbox{Ad}(C_t)(h) \in \mathfrak{n}$. From expression (\ref{expression2: horizontal lift}), it holds that:
\[
dC_t^{-1} = -R_{C_t^{-1}}V \diamond dZ_t.
\]
Hence, 
\[
dC_t = V^*(C_t) \diamond dZ_t.
\]
Then, $C_t$ is a solution of equation (\ref{decomp1}). By expression (\ref{expression2: horizontal lift}), 
\begin{eqnarray}
\label{expression3: horizontal lift}
(L_{C_t}dC_t^{-1})^*(\eta_t) = -(\mbox{Ad}(C_t)V)^*(\eta_t) \diamond dZ_t.
\end{eqnarray}
Finally, combining expressions (\ref{expression: horizontal lift}) and (\ref{expression3: horizontal lift}), it follows that:
\[
d\eta_t = (\mbox{Ad}(C_t)h)^*(\eta_t) \diamond dZ_t.
\]

\eop

\end{document}